\documentclass[reqno]{amsart}

\usepackage{amsmath}
\usepackage{amssymb}
\usepackage{bbm}
\usepackage{mathrsfs}

\makeatletter
\def\hsmash{\relax 
  \ifmmode\def\next{\mathpalette\mathhsm@sh}\else\let\next\makehsm@sh
  \fi\next}
\def\makehsm@sh#1{\setbox\z@\hbox{#1}\finhsm@sh}
\def\mathhsm@sh#1#2{\setbox\z@\hbox{$\m@th#1{#2}$}\finhsm@sh}
\def\finhsm@sh{\wd\z@\z@ \box\z@}
\makeatother

\newcommand{\br}{ }
\newcommand{\brr}{, }

\newtheorem{theo}{Theorem}[section]
\newtheorem{lem}[theo]{Lemma}
\newtheorem{prop}[theo]{Proposition}
\newtheorem{coro}[theo]{Corollary}
\newtheorem{obsc}[theo]{Conjecture}
\newtheorem{obst}[theo]{Theorem}

\theoremstyle{definition}
\newtheorem{ttt}[theo]{}
\newtheorem{defi}[theo]{Definition}
\newtheorem{algo}[theo]{Algorithm}
\newtheorem{str}[theo]{Strategy}

\newtheorem{rem}[theo]{Remark}
\newtheorem{rems}[theo]{Remarks}
\newtheorem{nota}[theo]{Notation}

\newtheorem{exs}[theo]{Examples}

\def\pmod#1{\nobreak\mkern8mu (\text{\rmfamily\upshape mod}\,\,#1)}
\renewcommand{\atop}[2]{\genfrac{}{}{0pt}{}{#1}{#2}}

\newcommand{\Pic}{\mathop{\text{\rm Pic}}\nolimits}
\newcommand{\Jac}{\mathop{\text{\rm Jac}}\nolimits}
\newcommand{\Gal}{\mathop{\text{\rm Gal}}\nolimits}
\newcommand{\Frob}{\mathop{\text{\rm Frob}}\nolimits}
\newcommand{\Kum}{\mathop{\text{\rm Kum}}\nolimits}
\newcommand{\End}{\mathop{\text{\rm End}}\nolimits}
\newcommand{\GL}{\mathop{\text{\rm GL}}\nolimits}
\newcommand{\GO}{\mathop{\text{\rm GO}}\nolimits}

\newcommand{\Tr}{\mathop{\text{\rm Tr}}\nolimits}

\newcommand{\im}{\mathop{\text{\rm im}}\nolimits}
\newcommand{\rk}{\mathop{\text{\rm rk}}\nolimits}
\newcommand{\spann}{\mathop{\text{\rm span}}\nolimits}
\newcommand{\disc}{\mathop{\text{\rm disc}}\nolimits}
\newcommand{\diag}{\mathop{\text{\rm diag}}\nolimits}

\newcommand{\bbC}{{\mathbbm C}}
\newcommand{\bbF}{{\mathbbm F}}
\newcommand{\bbN}{{\mathbbm N}}
\newcommand{\bbQ}{{\mathbbm Q}}
\newcommand{\bbZ}{{\mathbbm Z}}

\newcommand{\calJ}{{\mathscr{J}}}
\newcommand{\calO}{{\mathscr{O}}}

\newcommand{\frakA}{{\mathfrak{A}}}
\newcommand{\frakB}{{\mathfrak{B}}}
\newcommand{\frakE}{{\mathfrak{E}}}
\newcommand{\frakL}{{\mathfrak{L}}}
\newcommand{\frakX}{{\mathfrak{X}}}

\newcommand{\et}{{{\text{\rm {\'e}t}}}}
\newcommand{\tr}{{\text{\rm tr}}}
\newcommand{\ns}{{\text{\rm ns}}}

\newcommand{\GSp}{{\text{\rm GSp}}}
\newcommand{\USp}{{\text{\rm USp}}}

\newcommand{\Pb}{{\text{\bf P}}}
\newcommand{\Ab}{{\text{\bf A}}}

\newcommand{\Ib}{{\text{\bf I}}}

\newcounter{abc}

\newcounter{iii}
\newenvironment{iii}{\begin{list}{\rm \roman{iii}) }%
{\usecounter{iii} \leftmargin=0.0pt \labelsep=0.0pt %
\listparindent=0.0pt \labelwidth=0.0pt \parsep=\smallskipamount%
 \itemsep=0.0pt \topsep=0.0pt \partopsep=\smallskipamount}}{\end{list}}

\def\eop{\ifmmode\rule[-22pt]{0pt}{1pt}\ifinner\tag*{$\square$}\else\eqno{\square}\fi\else\vadjust{}\hfill$\square$\fi}

\begin{document}

\title[Examples of
$K3$~surfaces
with real multiplication]{Examples of
{\boldmath$K3$}~surfaces
with real multiplication}

\author{Andreas-Stephan Elsenhans}
\address{School of Mathematics and Statistics F07, University of Sydney, NSW 2006, Sydney, Australia, {\tt http://www.staff\!.\!uni-bayreuth.de/$\sim$bt270951/}}
\email{stephan@maths.usyd.edu.au}

\author{J\"org Jahnel}
\address{D\'epartement Mathematik, Universit\"at Siegen, Walter-Flex-Str.~3, D-57068 Siegen, Germany, {\tt http://www.uni-math.gwdg.de/jahnel}}
\email{jahnel@mathematik.uni-siegen.de}

%

\date{}
\dedicatory{}

\begin{abstract}
We construct explicit
$K3$
surfaces
over~$\bbQ$
having real~multiplication. Our~examples are of geometric Picard rank~16. The standard method for the computation of the Picard rank provably fails for the surfaces constructed.
\end{abstract}



\maketitle

\section{Introduction}

It~is well known that the endomorphism algebra of a general elliptic curve
$\frakX$
over
$\bbC$
is equal to
$\bbZ$,
while for certain exceptional curves the endomorphism algebra is~larger. There~are only countably many exceptions and these have complex multiplication~(CM). I.e.,
$\End(\frakX) \!\otimes_\bbZ\! \bbQ$
is an imaginary quadratic number~field.

There~is a rich theory about CM elliptic curves, cf.~\cite[Chapter~II]{Si} or~\cite[Chap\-ter~3]{Co}. We~will not go into details, but mention only a few facts that are relevant for what~follows. First~of all, the construction of CM elliptic curves in an analytic setting is very classical~\cite[17.~bis 23.~Abschnitt]{We}. The~situation becomes slightly more complicated, however, when explicit equations are asked~for.

For~$X$
an elliptic curve
over~$\bbQ$,
one says that it has complex multiplication if its base extension
$\frakX := X_\bbC$~has.
In~this situation, the occurrence of complex multiplication has striking consequences for the arithmetic
of~$X$.
For~example, on all general elliptic curves, the traces of the Frobenii
$\smash{\Frob_p \in \End(H^1_\et(X_{\overline\bbQ},\bbQ_l))}$
have the same statistic for
$p \to \infty$,
while, in the CM case, a different statistic~occurs.

To~be more precise, for a non-CM elliptic curve, the distribution of the normalized Frobenius traces
$\smash{\frac{\Tr\Frob_p}{2\sqrt{p}}}$
for
$p < N$
is supposed to converge, in the weak sense, to
$\frac2\pi\sqrt{1-t^2}\,dt$
when
$N \to \infty$.
Under~the additional assumption that
$X$
has at least one prime of multiplicative re\-duc\-tion, this behaviour has actually been established in 2010~\cite[Theorem~4.3]{HSBT}. On~the other hand, if
$X$
has CM by
$\bbQ(\sqrt{-d})$
then
$\Tr\Frob_p = 0$
for all primes that are inert 
in~$\bbQ(\sqrt{-d})$.

Further,~there are only nine imaginary quadratic number fields that may occur as the endomorphism field of a CM elliptic curve, defined
over~$\bbQ$,
namely those of class number~one.\medskip

The~whole theory generalizes to higher~dimensions. The~most obvious situation is certainly that of an abelian~surface. Here,~once again, the general case is that the endomorphism algebra is equal
to~$\bbZ$.

However, in contrast to the case of elliptic curves, there is more than one way for the endomorphism algebra to be~exceptional. For~instance, an abelian surface may have real multiplication (RM) \cite{Hu}. I.e.,~the endomorphism algebra may be an order in a totally real number field
$\supsetneqq\!\bbQ$.
Concerning~the possible statistics of the Frobenii on an abelian surface
over~$\bbQ$,
interesting investigations have been undertaken by F.~Fit\'e, K.\,S.~Kedlaya, V.~Rotger, and A.V.~Sutherland \cite{FKRS,KS}.

The~four authors present evidence for the existence of in fact 52 distinct types of abelian surfaces. A~theoretical explanation for this is as~follows. Associated~to every abelian
surface~$X$
over~$\bbQ$,
there is an algebraic group
$\smash{G \subset \End(H^1_\et(X_{\overline\bbQ},\bbQ_l)) \cong \GSp_4(\bbQ_l)}$
such that the image of the natural operation of
$\Gal(\overline\bbQ/\bbQ)$
on
$\smash{H^1_\et(X_{\overline\bbQ},\bbQ_l)}$
is Zariski dense
in~$G$.
The~reader might compare the Theorem of Tankeev and Zarhin [Theorem~\ref{TZ}], which gives an analogous statement for
$K3$~surfaces.
Corresponding to
$G$,
there is a compact subgroup
of~$\USp_4$.

On~the other hand, up to conjugation,
$\USp_4$
has exactly 55 compact subgroups that fulfill a number of necessary conditions \cite[Definition~3.1]{FKRS}. For~52 of them, an actual abelian surface~exists. Among~these, however, only 34 may be realized by an abelian surface
over~$\bbQ$.
The~others need larger base fields \cite[Theorem~4.3]{FKRS}. Furthermore,~the idea that the Frobenius elements
$\Frob_p$
are in fact equidistributed with respect to the Haar measure leads to hy\-po\-thet\-i\-cal distributions for the normalized Frobenius traces, which seem to agree with experimental observations.\medskip

From~the point of view of the classification of algebraic surfaces~\cite{Be}, abelian surfaces are not the only kind naturally generalizing elliptic curves to dimension two. Another~is provided by the so-called
$K3$~surfaces.
Indeed, elliptic curves may be characterized by the properties that they are of dimension one and have a trivial canonical~sheaf. On~the other hand, a surface with trivial canonical sheaf is either abelian
or~$K3$.

In~the case of an elliptic curve or abelian surface, the endomorphism field
$\End(\frakX) \!\otimes_\bbZ\! \bbQ$
is canonically isomorphic to the endomorphism field
$\End(H)$
of the associated Hodge structure
$H := H^1(\frakX, \bbQ)$.
This~may be just an equivalent reformulation, yet it allows to carry over the concepts of real and complex multiplication to more general~varieties. In~the particular situation of a
$K3$~surface,
the usual convention is to consider the field
$\End(T)$
of endomorphisms, in the category of Hodge structures [cf.~Section~\ref{sec_Hodgestruc}], of the transcendental part
$T \subset H^2(\frakX, \bbQ)$
of the second cohomology vector space [Example~\ref{ex_polarized}.ii)].\bigskip

\noindent
{\em B.~van Geemen's analytic approach.}\hskip3mm
B.~van Geemen showed that there exists a one-pa\-ram\-e\-ter family of
$K3$~surfaces
of Picard rank~16 that have real multiplication
by~$\smash{\bbQ(\sqrt{d})}$,
as soon as
$d$
is an odd integer that is a sum of two squares \cite[Example~3.4]{Ge}.

Van~Geemen's approach is analytic and does not lead to explicit equations. He~poses the problem to construct explicit examples in~\cite[Paragraph~3.1]{Ge}.
We~shall give van~Geemen's argument in a slightly more general form in an~appendix. In~fact, we will show that his method works for every
integer~$d$,
being even or odd, that is a sum of two squares. We~will also show that the four-dimensional part of the moduli stack of
$K3$~surfaces
he considered does not contain any surface having real multiplication
by~$\smash{\bbQ(\sqrt{d})}$,
when
$d$
is not a sum of two~squares.

\begin{ttt}
{\bf The results.}
In this note, we will present algorithms to efficiently test a
$K3$
surface~$X$
over~$\bbQ$
for real multiplication. Our~algorithms do not provide a proof, but only strong~evidence. Experiments~using them delivered two families of
$K3$~surfaces
of geometric Picard rank
$16$
and an isolated~example.

For~infinitely many members
$X^{(2,t)}$
of the first family, we will prove [Theorems~\ref{congr_Qw2_familie} and~\ref{RM_Qw2_familie}] that they have real multiplication
by~$\smash{\bbQ(\sqrt{2})}$.
To~our knowledge, these are the first explicit examples of
$K3$~surfaces
for which real multiplication is~proven.

The~members of the second family
$X^{(5,t)}$
are highly likely to have real multiplication
by~$\smash{\bbQ(\sqrt{5})}$,
while the isolated example
$X^{(13)}$
is strongly suspected to have real multiplication
by~$\smash{\bbQ(\sqrt{13})}$.\medskip

\noindent
{\em Our approach.}\hskip3mm
There is a theoretical algorithm to prove real (or complex) multiplication for a given
$K3$
surface under the assumption of the Hodge~conjecture. Cf.~the indications given in the proof of \cite[Theorem~6]{Ch1}. Its~main idea is to inspect the Hilbert scheme of
$X \times X$,
which is far from realistic to be done in~practice.

That~is why we decided to choose a different, more indirect approach. We~searched for surfaces having real multiplication through its arithmetic consequences. The~main idea behind our approach is that real multiplication
by~$\smash{\bbQ(\sqrt{d})}$
implies 
$\#X_p(\bbF_{\!p}) \equiv 1 \pmod p$,
for all
primes~$p$
that are inert under the field extension
$\smash{\bbQ(\sqrt{d})/\bbQ}$
[Corollary~\ref{RM_nonord}.i)]. This~result is in close analogy with the classical case of a CM elliptic curve and leads to an algorithm that is extremely selective, cf.~Section~\ref{sec_alg}.

From the surfaces found, we could guess the two families. For~the members of the first family, we will give a formal proof that
$\smash{\#X^{(2,t)}_p(\bbF_{\!p}) \equiv 1 \pmod p}$
is true for all primes
$p \equiv 3,5 \pmod 8$,
not just for those within the computational~range. In~order to do this, we analyze in detail one of the elliptic fibrations the surfaces
$X^{(2,t)}$
have~[Theorem~\ref{point_count}]. The~infinitely many congruences
$\smash{\#X^{(2,t)}_p(\bbF_{\!p}) \equiv 1 \pmod p}$
are enough to imply that the endomorphism field
$\End(T)$
is strictly larger than
$\bbQ$~[Lemma~\ref{congr_genuegt}].
For~infinitely many of the surfaces, actually
$\smash{\End(T) \cong \bbQ(\sqrt{2})}$
[Theorem~\ref{RM_Qw2_familie}].

In~the case of the family
$X^{(5,t)}$
and for the
surface~$X^{(13)}$,
we do not have a proof for the congruences on the point~count. The~experimental evidence is, however, overwhelming, cf.\ Remark~\ref{experiment_evidence}.
\end{ttt}

\begin{ttt}
{\bf An application. The analysis of van Luijk's method.}
Van Luijk's method is the standard method to determine the geometric Picard rank of a
$K3$~surface
over~$\bbQ$.
Its~fundamental idea is that, for every
prime~$p$
of good reduction, one has
$\smash{\rk \Pic X_{\overline\bbQ} \leq \rk \Pic X_{\overline\bbF_{\!p}}}$.
Further,~the method relies on the hope to find good primes such that
\begin{equation}
\label{eq_goodpr}
\smash{\rk \Pic X_{\overline\bbF_{\!p}} \leq \rk \Pic X_{\overline\bbQ} + 1} \, .
\end{equation}
To~see the method at work, the reader is advised to consult the original papers of R.~van Luijk~\cite{vL05,vL07} or some of the authors' previous articles~\cite{EJ08,EJ12,EJ13}. Further,~there is the remarkable work of N.~Elkies and A.~Kumar~\cite{EK}, in which they compute, among other data, the N\'eron-Severi ranks of all Hilbert-Blumenthal surfaces corresponding to the real quadratic fields of discriminants up
to~$100$.
Several~of them
are~$K3$.

Quite~recently, F.~Charles~\cite{Ch1} provided a theoretical analysis on the existence of primes fulfilling condition~(\ref{eq_goodpr}). The~result is that such primes always exist, unless
$X$
has real mul\-ti\-pli\-ca\-tion by a number
field~$E$
such that
$\smash{(22 - \rk\Pic X_{\overline\bbQ}) / [E\!:\!\bbQ]}$
is~odd. Thus,~our results provide explicit examples of
$K3$~surfaces
for which the method is bound to fail in its original~form.\smallskip

\looseness-1
Actually,~there is a more general version of van Luijk's method that applies to
$K3$~surfaces
having real multiplication, cf.~\cite[Proposition~18]{Ch1}. We~will make use of this in the proof of Theorem~\ref{RM_Qw2_familie}. It~works when the entire endomorphism field
$\End(T)$
is known, in particular, when
$\End(T) \cong \bbQ$.
Up~to now, no practical method has been found that would determine the geometric Picard rank of a
$K3$~surface
that has real multiplication, but for which this fact is not~known.
\end{ttt}

\section{Hodge structures}
\label{sec_Hodgestruc}

Recall the following definition, cf.~\cite[D\'efinition~2.1.10 and Proposition~2.1.9]{De71}.

\begin{defi}
i)
A (pure
$\bbQ$-)
{\em Hodge structure\/} of weight
$i$
is a finite dimensional
$\bbQ$-vector
space
$V$
together with a decomposition
$$V_\bbC := V \!\otimes_\bbQ\! \bbC = H^{0,i} \oplus H^{1,i-1} \oplus \ldots \oplus H^{i,0} \, ,$$
having the property that
$\overline{H^{m,n}} = H^{n,m}$
for every
$m,n \in \bbN_0$
such that
$m+n=i$.
A~{\em morphism\/}
$f\colon V \to V'$
{\em of (pure\/
$\bbQ$-)
Hodge structures\/} is a
$\bbQ$-linear
map such that
$f_\bbC\colon V_\bbC \to V'_\bbC$
respects the decompositions.\smallskip

\noindent
ii)
A~Hodge structure of
\mbox{weight~$2$}
is said to be {\em of
$K3$~type\/}
if
$\dim_\bbC H^{2,0} = 1$.
\end{defi}

\begin{rem}
Hodge structures of weight
$i$
form an abelian category~\cite[2.1.11]{De71}. Further,~this category is semisimple. I.e.,~every Hodge structure is a direct sum of primitive ones \cite[D\'efinition~2.1.4 and Proposition~2.1.9]{De71}.
\end{rem}

\begin{exs}
i)
Let~$\frakX$
be a smooth, projective variety
over~$\bbC$.
Then~$H^i(\frakX(\bbC),\bbQ)$
is in a natural way a pure
$\bbQ$-Hodge
structure of
weight~$i$.\smallskip

\noindent
ii)
In~$H^2(\frakX(\bbC),\bbQ)$,
the image of
$c_1\colon \Pic(\frakX) \!\otimes_\bbZ\! \bbQ \to H^2(\frakX(\bbC),\bbQ)$
defines a sub-Hodge
struc\-ture~$P$
such
that~$H^{0,2}_P = H^{2,0}_P = 0$,
which is called the {\em algebraic part\/} of
$H^2(\frakX(\bbC),\bbQ)$.
\end{exs}

\begin{defi}
i)
A {\em polarization\/} on a pure
$\bbQ$-Hodge
structure
$V$
of even weight is a nondegenerate symmetric bilinear form
$\langle.,.\rangle\colon V \times V \to \bbQ$
such that its
$\bbC$-bilinear
extension
$\langle.,.\rangle\colon V_\bbC \times V_\bbC \to \bbC$
satisfies the following two~conditions.

\begin{iii}
\item[ $\bullet$ ]
One has
$\langle x,y\rangle = 0$
for all
$x \in H^{m,n}$
and
$y \in H^{m',n'}$
such that
$m \neq n'$.
\item[ $\bullet$ ]
The inequality
$i^{m-n}\langle x,\overline{x}\rangle > 0$
is true for every
$0 \neq x \in H^{m,n}$.
\end{iii}

\noindent
ii)
A Hodge structure together with a polarization is called a {\em polarized Hodge structure.}
\end{defi}

\begin{exs}
\label{ex_polarized}
i)
If~$\frakX$
is a smooth, projective surface then
$H := H^2(\frakX(\bbC),\bbQ)$  
is a polarized pure Hodge structure, the polarization
$\langle.\,,.\rangle \colon H \times H \to \bbQ$
being given by the cup product, together with Poincar\'e duality.\smallskip

\noindent
ii)
The~algebraic part
$P \subseteq H$
and its orthogonal complement
$T = P^\perp$,
which is called the {\em transcendental part\/}
of~$H$,
are polarized Hodge structures, too.
If~$X$
is a
$K3$~surface
[Section~\ref{sec_K3}] then
$H$
and~$T$
are of
$K3$~type.
\end{exs}

\begin{ttt}
Yu.\,G.~Zarhin~\cite[Theorem~1.6.a) and Theorem~1.5.1]{Za83} proved that,
for~$T$
a polarized
\mbox{weight-$2$}
Hodge structure of
$K3$~type,
$E := \End(T)$
is either
$\bbQ$,
or a totally real field
$\supsetneqq \!\bbQ$,
or a CM~field.

Further,~every
$\varphi \in E$
operates as a self-adjoint~mapping. I.e.,
$$\langle \varphi(x), y \rangle = \langle x, \overline\varphi(y) \rangle \, ,$$
for
$\overline{\ \mathstrut}$
the identity map in the case that
$E$
is totally real and the complex conjugation in the case that
$E$
is a CM~field.

Observe that, in either case,
$T$
carries a structure of an
\mbox{$E$-vector}~space.
If~$E$
is totally real then one automatically has
$\dim_E T > 1$~\cite[Remark~1.5.3.c)]{Za83}.
\end{ttt}

\begin{defi}
Let~$T$
be a polarized
\mbox{weight-$2$}
Hodge structure of
$K3$~type.
If
$\End(T) \supsetneqq \bbQ$
is a totally real field then
$T$
is said to have {\em real multiplication}.
If~$\End(T)$
is CM then one speaks of {\em complex multiplication}.
\end{defi}

\section{Some background on
$K3$~surfaces}
\label{sec_K3}

\begin{ttt}
By~definition, a
$K3$~surface
is a simply connected, projective algebraic surface with trivial canonical~class.
\end{ttt}

\begin{exs}
Examples~include the classical Kummer surfaces, smooth space quartics and double covers
of~$\Pb^2\!$,
branched over a smooth sextic~curve. As~long as the singularities are isolated and rational, the minimal resolutions of singular quartics and double covers
of~$\Pb^2\!$,
branched over a singular sextic, are
$K3$~surfaces,
too.
In~this paper, we shall entirely work with the case of a double cover
of~$\Pb^2\!$,
branched over a singular~sextic.
\end{exs}

\begin{ttt}
The~property of being
$K3$
determines the Hodge~diamond. One~has
$H^1(\frakX, \bbQ) = 0$,
but
$H := H^2(\frakX, \bbQ)$
is non-trivial. It~is a pure
\mbox{weight-$2$}
Hodge structure of
dimension~$22$.
Further,
$\dim_\bbC H^{2,0}(\frakX) = \dim_\bbC H^{0,2}(\frakX) = 1$
and~$\dim_\bbC H^{1,1}(\frakX) = 20$.
The~Picard~group of a complex
$K3$~surface
is isomorphic
to~$\bbZ^n$,
where
$n$~may
range from
$1$
to~$20$.
\end{ttt}

\begin{defi}[{cf.~\cite[Paragraph~1.1]{Za83}}]
i)
Let~$\frakX$
be a complex
$K3$~surface
and
$T$
be the transcendental part of
$H^2(\frakX,\bbQ)$.
Then~$\frakX$
is said to have {\em real\/} or {\em complex multiplication\/} if
$T$~has.\smallskip

\noindent
ii)
A
$K3$~surface
$X$
over~$\bbQ$
is said to have {\em real\/} or {\em complex multiplication\/} if its base extension
$X_\bbC$~has.
\end{defi}

\begin{rems}
i)
The Kummer surface
$\Kum(\frakE_1 \!\times\! \frakE_2)$
attached to the product of two elliptic curves
$\frakE_1$
and
$\frakE_2$
has complex multiplication if one of the elliptic curves~has.\smallskip

\noindent
ii)
On~the other hand, a Kummer surface does {\em not\/} inherit the property of having real mul\-ti\-pli\-ca\-tion from the underlying abelian
surface~$\frakA$.
Indeed,~in this case,
$\smash{\sqrt{d} \in \bbQ(\sqrt{d})}$
operates on
$H^{1,0}(\frakA,\bbC)$
with
eigenvalues~$\smash{\pm\sqrt{d}}$.
Consequently,~$\smash{\bbQ(\sqrt{d})}$
operates on
$$\Lambda^2 H^{1,0}(\frakA,\bbC) = H^{2,0}(\frakA,\bbC) \hookrightarrow H_\bbC := H^2(\Kum\frakA,\bbC)$$
via multiplication by the norm, and the same is true for the whole
$T_\bbC \subset H_\bbC$.
\end{rems}

\begin{rem}
Motivated by the analysis of F.~Charles \cite{Ch1}, we are interested in
$K3$~surfaces
having real~multiplication and an
odd
\mbox{$E$-dimensional}~$T$.
The~simplest possible case is that
$\smash{E = \bbQ(\sqrt{d})}$
is real quadratic and
$\dim_E T = 3$,
i.e.~$\dim_\bbQ T = 6$.
\end{rem}

\begin{ttt}
{\em Frobenius eigenvalues.}\hskip3mm
For~varieties over finite fields, there is the
\mbox{$l$-adic}
cohomology theory~\cite{SGA5}.
If~$Y$
is a
$K3$~surface
over
$\bbF_{\!p}$
then
$\smash{\dim H^2_\et(Y_{\overline\bbF_{\!p}},\bbQ_l) = 22}$.
This~vector space is acted upon by
$\langle\Frob\rangle = \Gal(\overline\bbF_{\!p}/\bbF_{\!p})$.
The~22 eigenvalues are algebraic integers, independent of the choice of
$l \neq p$.
They are of absolute value
$p$
and
\mbox{$l$-adic}
units for every
$l \neq p$~\cite[Th\'eor\`eme~1.6]{De74}.

Concerning the
\mbox{$p$-adic}
nature of the Frobenius eigenvalues, there is the general result that the Newton polygon always runs above the Hodge polygon \cite{Ma}, cf.~\cite[Theorem~8.39]{BO}. A~variety over
$\bbF_{\!p}$
is called {\em ordinary\/} if the two polygons coincide~\cite[D\'e\-fi\-ni\-tion~IV.4.12]{IR}, cf.~\cite[pages~48f]{Il91}.

In~the particular case of a
$K3$~surface,
ordinarity is equivalent to the situation that the Frobenius eigenvalues are of
\mbox{$p$-adic}
valuations
$0,1,\ldots,1,2$.
On~the other hand, non-ordinarity implies that no Frobenius eigenvalue is a
\mbox{$p$-adic}~unit,
cf.~\cite[Paragraph~3.6]{Li}.
Therefore,~according to the Lefschetz trace formula~\cite[Expos\'e~XII, 6.3 and Exemple 7.3]{SGA5}, a
$K3$~surface
$Y$
over~$\bbF_{\!p}$
is ordinary if and only if
$\#Y(\bbF_{\!p}) \not\equiv 1 \pmod p$.
\end{ttt}

\section{Some arithmetic consequences of real multiplication}

Let~$X$
be a
$K3$~surface
over
$\bbQ$.
As~above, we~put
$$P := \im (c_1\colon \Pic(X_\bbC) \!\otimes_\bbZ\! \bbQ \hookrightarrow H^2(X(\bbC),\bbQ)) \, ,$$
$T := P^\perp$,
and
write~$E$
for the endomorphism algebra of the Hodge
structure~$T$.

Further,~let us choose a prime number
$l$
and turn to
$l$-adic
cohomology. This~essentially means to tensor
with~$\bbQ_l$,
as there is the canonical comparison isomorphism~\cite[Expos\'e~XI, Th\'eor\`eme~4.4.iii)]{SGA4}
$$H^2(X(\bbC),\bbQ) \!\otimes_\bbQ\! \bbQ_l \stackrel{\cong}{\longleftarrow} H^2_\et(X_{\overline\bbQ},\bbQ_l) \, .$$
An~important feature of the
\mbox{$l$-adic}
cohomology theory is that it is acted upon by the absolute Galois group of the base~field. I.e.,~there is a continuous representation
$$\varrho_l\colon \Gal(\overline\bbQ/\bbQ) \longrightarrow \GL(H^2_\et(X_{\overline\bbQ},\bbQ_l)) \, .$$
The image
of~$\varrho_l$
is an
$l$-adic
Lie group. Its~Zariski closure is an algebraic
group~$G_l$,
called the algebraic monodromy group associated
to~$\varrho_l$.

On~the other hand, there are the image
$P_l \subseteq H^2_\et(X_{\overline\bbQ},\bbQ_l)$
of
$\Pic(X_{\overline\bbQ}) \!\otimes_\bbZ\! \bbQ_l$
under the Chern map to
\mbox{$l$-adic}
cohomology
and its orthogonal
complement~$T_l$.
These~are compatible with the analogous constructions in Betti cohomology in that sense that
$P_l$
and
$T_l$
are mapped onto
$P \!\otimes_\bbQ\! \bbQ_l$
and
$T \!\otimes_\bbQ\! \bbQ_l$,
respectively, under the canonical comparison~isomorphism.

The~image
of~$\varrho_l$,
and hence the whole
of~$G_l$,
consists of endomorphisms
of~$H^2_\et(X_{\overline\bbQ},\bbQ_l)$
mapping
$P_l$
to~$P_l$.
Further,~these preserve orthogonality with respect to the
pairing~$\langle.\,,.\rangle$.
Thus,~the algebraic monodromy group
$G_l$
must map
$T_l$
into itself, as~well.

\begin{theo}[Tankeev, Zarhin]
\label{TZ}
The neutral component\/
$G_l^\circ$
of the algebraic monodromy group with respect to the Zariski topology is equal to the centralizer
of\/~$E$
in\/~$\GO(T_l, \langle.\,,.\rangle)$.
In~particular, the operation of\/
$E$
on\/
$\smash{T_l \subset H^2_\et(X_{\overline\bbQ},\bbQ_l)}$
commutes with that
of\/~$G_l^\circ$.\medskip

\noindent
{\bf Proof.}
{\em
This follows from the Mum\-ford-Tate conjecture, proven by S.\,G.~Tankeev \cite{Ta90,Ta95}, together with Yu.\,G.~Zarhin's explicit description of the Mumford-Tate group in the case of a
$K3$~surface~\cite[Theorem~2.2.1]{Za83}.
We~refer the reader to the original articles and to the discussion in~\cite[Section~2.2]{Ch1}.
}
\eop
\end{theo}

For~every
prime~$p$,
choose an absolute Frobenius element
$\smash{\Frob_p \in \Gal(\overline\bbQ/\bbQ)}$.
If~$p \neq l$
is a prime, at which
$X$
has good reduction then,~by virtue of the smooth base change theorem \cite[Exp.~XVI, Corollaire~2.5]{SGA4}, there is a canonical isomorphism
$$H^2_\et(X_{\overline\bbQ},\bbQ_l) \cong H^2_\et((X_p)_{\overline\bbF_{\!p}},\bbQ_l) \, .$$
Here,~the vector space on the right hand side is naturally acted upon by
$\Gal(\overline\bbF_{\!p}/\bbF_{\!p})$
and the operation of
$\smash{\Frob_p \in \Gal(\overline\bbQ/\bbQ)}$
on the left hand side is compatible with that of
$\smash{\Frob \in \Gal(\overline\bbF_{\!p}/\bbF_{\!p})}$
on the~right.

\begin{coro}
\label{commut}
There is a positive integer\/
$f$
such that, for every pair\/
$(p, l)$
of prime numbers, the operation of\/
$(\Frob_p)^f$
on\/~$T_l$
commutes with that
of\/~$E$.\medskip

\noindent
{\bf Proof.}
{\em
By~definition,
$\varrho_l(\Frob_p) \in G_l$.
Hence,~for
$f := \#(G_l/G_l^\circ)$,
we have
$\varrho_l((\Frob_p)^f) \in G_l^\circ$.
Further,~the groups
$G_l/G_l^\circ$
are canonically isomorphic to each other, for the various values
of~$l$,
as was proven by M.~Larsen and R.~Pink \cite[Propo\-si\-tion~6.14]{LP}.
}
\eop
\end{coro}

\begin{nota}
\label{poldef}
i)
For~every
prime~$p$,
choose
$l \neq p$
and denote by
$\smash{\chi^T_{p^n}}$
the characteristic polynomial of
$(\Frob_p)^n$
on the transcendental
part~$T_l$.
This~has coefficients
in~$\bbQ$
and is in\-de\-pen\-dent
of~$l$,
whether
$X$
has good reduction
at~$p$
\cite[Th\'eor\`eme~1.6]{De74} or not~\cite[Theorem~3.1]{Oc}.
One has
$\smash{\deg\chi^T_{p^n} = 22 - \rk\Pic X_{\overline\bbQ}}$.\smallskip

\noindent
ii)
We factorize
$\smash{\chi^T_{p^n} \in \bbQ[Z]}$
in the form
$$\chi^T_{p^n}(Z) = \chi^\tr_{p^n}(Z) \cdot \prod_{k,i} (Z - \zeta_k^i)^{e_{k,i}} ,$$\vskip-2mm\noindent
for
$\zeta_k := \exp(2\pi i/k)$,
$e_{k,i} \geq 0$,
and such that
$\smash{\chi^\tr_{p^n} \in \bbQ[Z]}$
does not have any roots of the form
$p^n$
times a root of~unity.
\end{nota}

\begin{rem}
If~$p$
is a good prime then, according to the Tate conjecture,
$\smash{\chi^\tr_{p^n}}$
is the characteristic polynomial of
$\Frob^n$
on the transcendental part of
$\smash{H^2_\et(X_{\overline\bbF_{\!p}}\!,\bbQ_l)}$.
In particular,
$\smash{\deg\chi^\tr_{p^n} = 22 - \rk\Pic X_{\overline\bbF_{\!p}}}$.
Further,
$\smash{\chi^\tr_{p^n} = \chi^T_{p^n}}$
if and only if
$\smash{\rk\Pic X_{\overline\bbF_{\!p}} = \rk\Pic X_{\overline\bbQ}}$.
\end{rem}

For~the remainder of this section, we assume that
$E \supseteq \bbQ(\sqrt{d})$,
for
$d \neq 1$
a square-free~integer. I.e.,~that
$X$
has real or complex multiplication by a number
field~$E$
that
contains~$\bbQ(\sqrt{d})$.
Further,~we shall use the symbol
$f$
for an arbitrary positive integer such that the operation of
$(\Frob_p)^f$
on~$T_l$
commutes with that of
$\bbQ(\sqrt{d})$
[Corollary~\ref{commut}].

\begin{prop}
\label{factoriz_l_ad}
Let\/~$p$
be a prime number and\/
$l$
be a prime that is ramified or inert
in\/~$\smash{\bbQ(\sqrt{d})}$.
Then~the polynomial\/
$\smash{\chi^\tr_{p^f} \in \bbQ[Z]}$
splits~as
$$\chi^\tr_{p^f} = g_lg_l^\sigma ,$$
for\/
$g_l \in \bbQ_l(\sqrt{d})[Z]$
and\/
$\smash{\sigma\colon \bbQ_l(\sqrt{d}) \to \bbQ_l(\sqrt{d})}$
the~conjugation.\medskip

\noindent
{\bf Proof.}
{\em
The~assumption ensures that
$\smash{\bbQ_l(\sqrt{d})}$
is a quadratic extension~field.
Further,~$T_l$
is a
$\smash{\bbQ_l(\sqrt{d})}$-vector
space and, by Corollary~\ref{commut},
$\smash{\varrho_l((\Frob_p)^f)}$
commutes with the operation of
$\sqrt{d} \in E$.
In~other words,
$\smash{\varrho_l((\Frob_p)^f)}$
is a
$\bbQ_l(\sqrt{d})$-linear~map.

For~the corresponding char\-ac\-ter\-is\-tic poly\-no\-mial
$\smash{c_l \in \bbQ_l(\sqrt{d})[Z]}$,
we have
$\smash{\chi^T_{p^f} = c_lc_l^\sigma}$.
The~assertion immediately follows from~this.
}
\eop
\end{prop}

\begin{lem}
\label{split}
Let\/~$K$
be any field,
$K(\sqrt{d})/K$
a quadratic field extension, and\/
$h \in K[Z]$
an irreducible~polynomial.
Then\/~$h$
splits
over\/~$K(\sqrt{d})$
if and only if\/
$K(\sqrt{d}) \subseteq K[Z]/(h)$.\medskip

\noindent
{\bf Proof.}
{\em
Suppose first that
$K(\sqrt{d}) \subseteq K[Z]/(h)$
and let
$z_0 \in K[Z]/(h)$
be a root
of~$h$.
Then
$K[Z]/(h) \cong K(z_0)$
and
$\smash{[K(z_0) \!:\! K(\sqrt{d})] = \frac{\deg h}2}$.
Therefore, the minimal polynomial
of~$z_0$
over~$K(\sqrt{d})$
is of degree~$\smash{\frac{\deg h}2}$
and a factor
of~$h$.

On the other hand, assume that
$h$
splits
over~$K(\sqrt{d})$
and
write~$h = gg^\sigma$.
Then~the extension fields
$K[Z]/(h)$
and
$K(\sqrt{d})[Z]/(g)$
both contain a zero
of~$g$
and have the same degree
over~$K$.
Hence,~they must be isomorphic to each~other.
}
\eop
\end{lem}

\begin{nota}
For~$e \in \bbN$
and a normalized polynomial
$h \in \bbQ[Z]$,
we will write
$h^{(e)}$
to denote the normalized polynomial of the same degree
as~$h$
that has the zeroes
$\smash{x_1^e, \ldots, x_r^e}$,
for
$x_1, \ldots, x_r$
the zeroes
of~$h$,
taken with~multiplicities.
\end{nota}

\begin{rems}
i)
For an irreducible polynomial
$h \in \bbQ[Z]$,
the
polynomial~$h^{(e)}$
must not factor, except as the power of an irreducible~polynomial. In~fact,
$\Gal(\overline\bbQ/\bbQ)$
permutes the roots
$x_1,\ldots,x_r$
of~$h$
transitively. Therefore,~it does the same to
$\smash{x_1^e,\ldots,x_r^e}$.\smallskip

\noindent
ii)
If~$h \in \bbQ[Z]$
is irreducible of
degree~$r$
and
$h(\zeta Z) \neq \zeta^r h(Z)$
for every
$e$-th
root of unity
$\zeta$
then
$h^{(e)}$
is~irreducible.
\end{rems}

\begin{theo}
\label{factoriz_rat}
Let\/~$p$
be a prime of good reduction of the\/
$K3$~surface\/~$X$
over\/~$\bbQ$,
having real or complex multiplication by a
field\/~$E$
containing the quadratic number
field\/~$\smash{\bbQ(\sqrt{d})}$.
Then~at least one of the following two statements is~true.

\begin{iii}
\item
The polynomial\/
$\chi^\tr_p \in \bbQ[Z]$
splits in the~form
$$\chi^\tr_p = gg^\sigma,$$
for\/
$\smash{g \in \bbQ(\sqrt{d})[Z]}$
and\/
$\smash{\sigma\colon \bbQ(\sqrt{d}) \to \bbQ(\sqrt{d})}$
the~conjugation.
\item
The polynomial\/
$\smash{\chi^\tr_{p^f}}$
is a square
in\/~$\bbQ[Z]$.
\end{iii}\medskip

\noindent
{\bf Proof.}
{\em
According~to \cite[Theorem~1.4.1]{Za93},
$\smash{\chi^\tr_p = h^k}$,
for an irreducible polynomial
$h \in \bbQ[Z]$
and
$k \in \bbN$.
Write
$\smash{h^{(f)} = \underline{h}^{k'}}$,
for
$\underline{h}$
an irreducible~polynomial.
Then~$\smash{\chi^\tr_{p^f} = \underline{h}^{kk'}}$.

If~one of the integers
$k$
and
$k'$
is even then assertion~ii) is~true. Thus,~assume from now on that
$k$
and
$k'$
are both~odd.
By~Proposition~\ref{factoriz_l_ad},
$\smash{\underline{h}^{kk'} = \chi^\tr_{p^f}}$
splits into two factors conjugate over
$\smash{\bbQ_l(\sqrt{d})}$,
for every
$l$
that is not split
in~$\smash{\bbQ(\sqrt{d})}$.
As~$kk'$
is odd, the same is true
for~$\underline{h}$.

In~particular, for every prime
$\frakL$
lying
above~$(l)$
in the field
$\bbQ[Z]/(\underline{h})$,
one has that
$f(\frakL|(l))$
is even for
$l$
inert
in~$\smash{\bbQ(\sqrt{d})}$
and that
$e(\frakL|(l))$
is even for
$l$
ramified
in~$\smash{\bbQ(\sqrt{d})}$.
\mbox{\cite[Chapter~VII, Proposition~13.9]{Ne}} implies that
$\smash{\bbQ(\sqrt{d}) \subseteq \bbQ[Z]/(\underline{h})}$.

Let~now
$\smash{x_0 \in \overline\bbQ}$
be an element having minimal
polynomial~$h$.
Then~$\smash{\bbQ(x_0^f) \cong \bbQ[Z]/(\underline{h})}$.
Altogether,~$\smash{\bbQ(\sqrt{d}) \subseteq \bbQ(x_0^f) \subseteq \bbQ(x_0)}$.
But,~according to Lemma~\ref{split}, this is equivalent to
$h$
being reducible over
$\smash{\bbQ(\sqrt{d})}$.
It~must split into two conjugate~factors.
}
\eop
\end{theo}

\begin{rems}
\label{block}
i)
Let~$h$
be an irreducible polynomial such that
$\smash{\chi^\tr_p = h^k}$
and consider
$\Gal(h)$
as a permutation group on the roots
of~$h$.
As~such, it has an obvious block structure \mbox{\cite[Section~1.5]{DM}}
$\frakB := \big\{ \{z, \overline{z}\} \mid h(z) = 0 \big\}$
into blocks of size~two.
Indeed,~$h$
is a real polynomial without real roots and every root is of absolute
value~$p$.
Thus,~$\smash{\overline{z} = \frac{p^2}z}$
and so the pairs are respected by the operation of the Galois~group.\smallskip

\noindent
ii)
Assume~that
$k$
is~odd and that
$d > 0$.
We~claim that this causes a second block~structure.
To~show this, let us suppose first that variant~i) of Theorem~\ref{factoriz_rat} is~true. Then~there is the block structure
$\frakB' := \big\{ \{z \mid g(z) = 0\}, \{z \mid g^\sigma(z) = 0\} \big\}$
into two blocks of size
$\smash{\frac{\deg h}2}$.
As~$g$
and~$g^\sigma$
are real polynomials, the blocks in
$\frakB'$
are non-minimal. Each~is a union of some of the blocks
in~$\frakB$.

If~option~ii) of Theorem~\ref{factoriz_rat} happens to be true then there is a block structure
$\frakB''$,
the blocks in which are formed by the roots
of~$h$
having their
\mbox{$f$-th}~power
in~common. The~mutual refinement of
$\frakB''$
and~$\frakB$
is the trivial block structure into blocks of size~one.
As~$k$
is assumed odd, the blocks
in~$\frakB''$
are of even~size. Thus,~the block structure generated by
$\frakB''$
and~$\frakB$
consists of blocks of a size that is a multiple
of~$4$.
\end{rems}

\begin{coro}
Suppose that\/
$d > 0$.
Then,~for every good
prime\/~$p$,
$\deg \chi^\tr_p$
is divisible
by\/~$4$.\medskip

\noindent
{\bf Proof.}
{\em
Write~$\smash{\chi^\tr_p = h^k}$.
As~seen in Remark~\ref{block}.i),
$\deg h$
is even, which implies the claim as long as
$k$
is~even.
When~$k$
is odd, the observations made in Remark~\ref{block}.ii) show in both cases that
$\deg h$
must be divisible
by~$4$.
}
\eop
\end{coro}

\begin{coro}
\label{2mod4}
Suppose that\/
$d > 0$.
Then,~for every good prime\/
$p \geq 3$,
we have
$$\rk \Pic((X_p)_{\overline\bbF_{\!p}}) \equiv 2 \pmod 4 \, .$$

\noindent\looseness-1
{\bf Proof.}
{\em
The~Tate conjecture is known to be true for
$K3$~surfaces
in
characteristic~$\geq \!3$,
cf.\ \mbox{\cite[Theorem~1]{Pe}}, \cite[Corollary~2]{Ch2}, and~\cite{LMS}. Further,~the characteristic polynomial of
$\Frob$
on
$\smash{H^2_\et((X_p)_{\overline\bbF_{\!p}}\!, \bbQ_l)}$
has exactly
$22 - \deg \chi^\tr_p$
zeroes of the form
$p$
times a root of~unity.
}
\eop
\end{coro}

\begin{coro}
\label{RM_nonord}
Suppose that\/
$d > 0$
and
let\/~$p \geq 3$
be a good prime number that is inert
in\/~$E = \bbQ(\sqrt{d})$.

\begin{iii}
\item
Then\/~$X_p$
is non-ordinary.
\item
Suppose~that\/
$\dim_E T \leq 3$.
Then,~either\/
$\smash{\rk \Pic((X_p)_{\overline\bbF_{\!p}}) = 22}$
or\/
$\chi^\tr_{p^f}$
is the square of an irreducible quadratic~polynomial.
\end{iii}\medskip

\looseness-1
\noindent
{\bf Proof.}
{\em
i)
$X_p$
being ordinary would mean that
$\smash{\chi^\tr_{p^f}}$
has exactly one zero that is a
\mbox{$p$-adic}~unit.
By~Theorem~\ref{factoriz_rat}, in any case, we can say that there is a factorization
$\smash{\chi^\tr_{p^f} = \underline{g}\underline{g}^\sigma}\!$,
for some
$\underline{g} \in \calO_E[Z]$.
Assume without restriction that the zero being a
\mbox{$p$-adic}~unit
is a root
of~$\underline{g}^\sigma$.
Then, for the coefficients of the~polynomial
$$\underline{g}(Z) = Z^n + a_{n-1}Z^{n-1} + \ldots + a_0 \, ,$$
one has that
$\nu_p(a_j) > 0$,
for
every~$j$.
But,~$p$
is inert, hence the same is true
for~$\underline{g}^\sigma$.
In par\-tic\-u\-lar,
$\nu_p(a_{n-1}^\sigma) > 0$.
This~shows that it is impossible for
$\underline{g}^\sigma$
to have exactly one root that is a
\mbox{$p$-adic}~unit.\smallskip

\noindent
ii)
The assumption
$\dim_E T \leq 3$
means
$\smash{\rk\Pic X_{\overline\bbQ} \geq 16}$.
Then, even more,
$\smash{\rk \Pic((X_p)_{\overline\bbF_{\!p}}) \geq 16}$.
From~Corollary~\ref{2mod4}, we know that either
$\smash{\rk \Pic((X_p)_{\overline\bbF_{\!p}}) = 18}$
or~$\smash{\rk \Pic((X_p)_{\overline\bbF_{\!p}}) = 22}$.
As~in the latter case the proof is complete, let us suppose that the rank
is~$18$.

Then~$\deg \chi^\tr_p = \deg \chi^\tr_{p^f} = 4$.
Theorem~\ref{factoriz_rat} gives us two options. Option~ii) is that
$\chi^\tr_{p^f} = g^2$
is the square of a quadratic polynomial
$g \in \bbQ[Z]$.
Since~its roots are non-reals, this polynomial must be~irreducible.

Otherwise,~according to option~i), there is a factorization
$\smash{\chi^\tr_{p} = gg^\sigma}$\!,
for
some~$g \in E[Z]$.
Write~$g(Z) = Z^2 + aZ \pm p^2 = (Z-x_1)(Z-x_2)$.
Then
$$\nu_p(x_1) + \nu_p(x_2) = \nu_p(x_1x_2) = \nu_p(\pm p^2) = 2$$
and
$\min(\nu_p(x_1), \nu_p(x_2)) \leq \nu_p(x_1 + x_2) = \nu_p(-a)$.
If~$\nu_p(x_1) \neq \nu_p(x_2)$
then equality is~true.

Further,~$\smash{a \in \bbQ(\sqrt{d})}$
implies that
$\nu_p(-a)$
is an integer and it is well-known that
$\nu_p(x_i) \geq 0$.
Thus,~there are only two~cases. We~will show that they are both~contradictory.

If~$\nu_p(x_1) = \nu_p(x_2) = 1$
then, as
$p$
is inert, the same is true for
$x_1^\sigma$
and~$x_2^\sigma$.
Hence,~the four quotients
$x_1/p$,
$x_2/p$,
$x_1^\sigma/p$,
and~$x_2^\sigma/p$
are
$p$-adic
units. On~the other hand, the eigenvalues of
$\Frob_p$
on
$l'$-adic
cohomology
are known to be
$l'$-adic
units for every
$l' \neq p$.
Hence,~$x_1/p$,
$x_2/p$,
$x_1^\sigma/p$,
and~$x_2^\sigma/p$
are
$l$-adic
units actually for all
primes~$l$.
Consequently,~they must be roots of~unity. This,~however, is a contradiction to the definition
of~$\smash{\chi^\tr_p}$,
given in~\ref{poldef}.ii).

On~the other hand, if, without restriction,
$\nu_p(x_1) = 0$
and
$\nu_p(x_2) = 2$
then
$\nu_p(x_1^\sigma) = 0$,~too.
This~is a contradiction to the general fact that the Newton polygon always runs above the Hodge~polygon.
}
\eop
\end{coro}

\begin{coro}
\label{Klein_four}
Suppose~that\/
$\dim_E T \leq 3$.
If\/~$\chi^\tr_{p^f}$
is the square of a quadratic poly\-nomial, but\/
$\chi^\tr_p$
is not, then\/
$\Gal(\chi^\tr_p) \cong \bbZ/2\bbZ \times \bbZ/2\bbZ$.\medskip

\noindent
{\bf Proof.}
{\em
The~assumption implies that
$\smash{\chi^\tr_p = h}$
is irreducible of degree~four.
Further,~$\Gal(h)$
has two different block structures, both into blocks of size~two. The~only transitive subgroup
of~$S_4$
having this property is the Klein four~group.
}
\eop
\end{coro}

\section{Efficient algorithms to test a
$K3$~surface
for real multiplication}
\label{sec_alg}

\paragraph*{{\em Generalities.}}
Recall~that a
$K3$~surface
$Y$
over a finite
field~$\bbF_{\!p}$
is ordinary if and only if
$\#Y(\bbF_{\!p}) \not\equiv 1 \pmod p$.
In~particular, non-ordinarity may be tested by counting points only
over~$\bbF_{\!p}$.

For
$K3$~surfaces
with real multiplication
by~$\smash{\bbQ(\sqrt{d})}$,
we expect non-ordinary reduction at approximately half the~primes. On~the other hand, consider a general
$K3$~surface
$X$
over~$\bbQ$
of a certain geometric Picard~rank. I.e.,~assume that
$\End(T) \cong \bbQ$.
Then~Theorem~\ref{TZ} implies that the Frobenii
$\varrho_l(\sigma\Frob_p\sigma^{-1})$,
for
$p$
running through the primes and
$\sigma$
through
$\Gal(\overline\bbQ/\bbQ)$,
are Zariski dense in
$\GO(T_l, \langle.\,,.\rangle)$.
In~particular, the values
$\smash{\frac{\#X_p(\bbF_{\!p}) - p^2 - 1}p = \frac1p\!\Tr\Frob_p}$
are Zariski dense
in~$\Ab^1$.
In~a way similar to \cite{KS}, one may hope that
$\smash{(\frac1p\!\Tr\Frob_p \bmod 1)}$
is equidistributed
in~$[0,1]$.

Thus,~somewhat naively, we expect that a general
$K3$~surface
$X$
over~$\bbQ$
has non-ordinary reduction
at~$p$
with a probability
of~$\smash{\frac1p}$.
The~number of primes
$\leq \!N$,
at which the reduction is non-ordinary, should be of the order of
$\log\log N$.

This~suggests to generate a huge sample of
$K3$~surfaces
over~$\bbQ$,
each having geometric Picard rank
$\geq \!16$,
and to execute the following statistical algorithm on all of~them.

\begin{algo}[Testing a
$K3$~surface
for real multiplication--statistical version]
\label{alg_stat}
\leavevmode
\begin{iii}
\item
Let~$p$
run over all primes
$p \equiv 1 \pmod 4$
between
$40$
and~$300$.
For~each~$p$,
count the number
$\#X_p(\bbF_{\!p})$
of
$\bbF_{\!p}$-rational
points on the reduction
of~$X$
modulo~$p$.
If~$\#X_p(\bbF_{\!p}) \equiv 1 \pmod p$
for not more than five primes then terminate immediately.
\item
Put~$p_0$
to be the smallest good and ordinary prime
for~$X$.
\item
Determine the characteristic polynomial of
$\Frob$
on
$\smash{H^2_\et((X_{p_0})_{\overline\bbF_{\!p_0}}\!, \bbQ_l)}$.
For this, use the strategy described in~\cite[Examples 27 and~28]{EJ10}. 
Factorize~the polynomial obtained to calculate the polynomial
$\smash{\chi_{p_0}^\tr}$.
If~$\smash{\deg \chi_{p_0}^\tr \neq 4}$
then~terminate.

Test~whether
$\smash{\chi_{p_0}^\tr}$
is the square of a quadratic polynomial. In~this case, raise
$p_0$
to the next good and ordinary prime and iterate this step.

Otherwise,~determine the Galois group
$\smash{\Gal(\chi_{p_0}^\tr)}$.
If~$\smash{\Gal(\chi_{p_0}^\tr)}$
is isomorphic to the Klein four group then raise
$p_0$
to the next good and ordinary prime and iterate this step.
\item
Now,~$\smash{\chi_{p_0}^\tr}$
is irreducible of degree~four. Determine the quadratic subfields of the splitting field
of~$\smash{\chi_{p_0}^\tr}$.
Only~one real quadratic field may~occur.
Put~$d$
to be the corresponding~radicand.
\item
Let~$p$
run over all primes
$< \!300$
that are inert
in~$\bbQ(\sqrt{d})$,
starting from the~lowest. If
$\#X_p(\bbF_{\!p}) \not\equiv 1 \pmod p$
for one these then~terminate.
\item
Output~a message saying that
$X$
is highly likely to have real or complex multiplication by a field containing~$\smash{\bbQ(\sqrt{d})}$.
\end{iii}
\end{algo}

\begin{rems}
\label{rems_alg_stat}
i)
Algorithm~\ref{alg_stat} does not give false negatives due to bad~reduction. Cf.\ Lemma \ref{congr_bad_red}, below.\smallskip

\noindent
Nevertheless,~the algorithm is only statistically~correct. It~is possible, in principle, that a
$K3$~surface
with real multiplication is thrown away in step~i). However,~in the case that
$\End(T) = \bbQ(\sqrt{d})$,
this may occur only if not more than five of the primes used in the algorithm are inert
in~$\bbQ(\sqrt{d})$.
The~smallest discriminant, for which this happens,
is~$d = 8493$.\smallskip

\noindent
ii)
On~the other hand, Algorithm~\ref{alg_stat} is extremely efficient. The point is that, for the lion's share of the surfaces, it terminates directly after step~i). In~fact, according to the inclusion-exclusion principle \cite[formula~(2.1.3)]{Ha}, the likelihood that a surface with
$\End(T) \cong \bbQ$
survives step~i) should~be
$$\sum_{r=6}^{\#S} \sum_{\atop{R \subset S}{\#\!R = r}} (-1)^{r-6} \Big( \atop{r}6 \Big) \frac1{\prod\limits_{p \in R}\!\!p} \approx 2.66 \!\cdot\! 10^{-8} \, ,$$
for~$S := \{p \mid p {\rm ~prime,} ~40<p<300, ~p \equiv 1 \pmod 4\}$. Thus,~only for a negligible percentage of the surfaces, the more time-consuming steps ii)-v) have to be carried~out.\smallskip

\noindent
This shows, in particular, that step~i) is the only time-critical~one. An~efficient algorithm for point counting over relatively small prime fields is asked~for.\smallskip

\noindent
iii)
In~our samples, step~iii) involves to count, in~addition, the points
on~$\smash{X_{p_0}}$
that are defined over
$\smash{\bbF_{\!p_0^2}}$
and, possibly,
over~$\smash{\bbF_{\!p_0^3}}$,
but not over larger~fields. The~reason for this is that 16 generators of the cohomology vector space are explicitly known, including the Galois operation on~them. Thus,~only a degree six factor of the desired polynomial of
degree~$22$
needs to be~computed.\smallskip

\noindent
iv)
The second part of step~iii) has the potential to create an infinite~loop. But~this never happened for any of the surfaces we tested. Whenever step~i) suggested real multiplication, after a few trials we found a
prime~$p_0$
such that
$\smash{\deg \chi_{p_0}^\tr}$
was irreducible of degree four and had the cyclic group of order four or the dihedral group of order eight as its Galois~group.\smallskip

\noindent
v)
In~step~iv), the polynomial
$\smash{\chi_{p_0^f}^\tr}$
is certainly irreducible although the value
of~$f$
is not known to~us. This~is simply the assertion of Corollary~\ref{Klein_four}. As~a consequence of this, Theorem~\ref{factoriz_rat} shows that
$\smash{\chi_{p_0}^\tr}$
must split over the RM~field.\smallskip

\noindent
vi)
The~reason for restricting in step~i) to primes congruent to
$1$
modulo~$4$
is a practical~one. Otherwise,~too many surfaces are found showing the pattern that
$\#X_p(\bbF_{\!p}) \equiv 1 \pmod p$
for every prime
$p \equiv 3 \pmod 4$.
These~primes are inert under
$\bbQ(\sqrt{-1})/\bbQ$,
but not under any real quadratic field~extension. We~do not exactly understand why our samples contained many more such surfaces than those we were looking~for.\smallskip

\noindent
On~the other hand, for small
primes~$p$,
it happens too often that
$\#X_p(\bbF_{\!p}) \equiv 1 \pmod p$,
independently of whether or not
$X$
has real~multiplication. As~this would slow down the algorithm, we incorporated the restriction to primes
$p > 40$.\smallskip

\noindent
vii)
The likelihood that a random surface would survive step~v) is
\begin{equation}
\label{eq_prob}
\prod_{\atop{p {\rm \,inert\,in\,} \bbQ(\sqrt{d}),}{p<300}} \hskip-5.5mm 1/p \, ,
\end{equation}
which is less than
$10^{-60}$
for small values
of~$d$.
Thus,~we do not expect any false positives to be given by Algorithm~\ref{alg_stat}.
\end{rems}

When testing surfaces for real multiplication by a particular
field~$\smash{\bbQ(\sqrt{d})}$,
the following mod\-i\-fi\-ca\-tion of Algorithm~\ref{alg_stat} may be~used.

\begin{algo}[\mbox{Testing a
$K3$~surface
for real multiplication--deterministic version}]
\label{alg_det}
\leavevmode
\begin{iii}
\item[o) ]
This~algorithm assumes that, in an initialization step, the primes
$p < 300$
have been listed that are inert
in~$\smash{\bbQ(\sqrt{d})}$.
\item
Let~$p$
run over the~list.
For~each~$p$,
count the numbers
$\#X_p(\bbF_{\!p})$
of
$\bbF_{\!p}$-rational
points on the reduction
$X_p$.
If~one of them is not congruent
to~$1$
modulo~$p$
then terminate~immediately.
\item
Let~$p$
run over all good primes
$< \!100$,
starting from the~lowest.

For~each prime, calculate the polynomial
$\smash{\chi_p^\tr}$,
as in Algorithm~\ref{alg_stat}.iii).
If~$\smash{\deg \chi_p^\tr \neq 0}$
or~$4$
then terminate.
If~$\smash{\deg \chi_p^\tr = 4}$
then test whether
$\smash{\chi_p^\tr}$
is the square of a quadratic polynomial. If~this is the case then go to the next~prime.

Factor
$\smash{\chi_p^\tr}$
over~$\smash{\bbQ(\sqrt{d})}$
and determine the Galois group
$\smash{\Gal(\chi_p^\tr)}$.
If~neither
$\smash{\chi_p^\tr}$
splits over
$\smash{\bbQ(\sqrt{d})}$
nor
$\smash{\Gal(\chi_p^\tr) \cong \bbZ/2\bbZ \times \bbZ/2\bbZ}$
then~terminate.
Otherwise,~go to the next~prime.
\item
Output~a message saying that
$X$
is highly likely to have real or complex multiplication by a field containing~$\smash{\bbQ(\sqrt{d})}$.
\end{iii}
\end{algo}

\begin{rems}
i)
Algorithm~\ref{alg_det} does not give false~negatives. Bad~reduction does not cause any problem, due to Lemma~\ref{congr_bad_red}.\smallskip

\noindent
ii)
The~likelihood that a general
$K3$~surface
survives step~i) is again given by formula~(\ref{eq_prob})~above. In~the cases
$d = 2$,
$5$,
$13$,
and
$17$,
where we actually run the algorithm, the values of the product are approximately
$3.26 \!\cdot\! 10^{-64}$,
$2.69 \!\cdot\! 10^{-63}$,
$4.07 \!\cdot\! 10^{-61}$,
and
$1.30 \!\cdot\! 10^{-63}$.
In~accordance with this, no statistical outliers showed up in step~ii).
\end{rems}

\begin{lem}
\label{congr_bad_red}
Let\/~$X$
be a double cover
of\/~$\Pb^2_\bbQ$,
branched over the union of six~lines. Suppose~there is a quadratic number field\/
$\smash{\bbQ(\sqrt{d})}$
such that\/
$\#X_q(\bbF_{\!q}) \equiv 1 \pmod q$
for every good prime\/
$q$
that is inert
in\/~$\smash{\bbQ(\sqrt{d})}$.\smallskip

\noindent
Then\/~$\#X_p(\bbF_{\!p}) \equiv 1 \pmod p$,
too, for every bad
prime\/~$p$
that is~inert.\medskip

\noindent
{\bf Proof.}
{\em
If~at least two of the six lines coincide
modulo~$p$
then
$X_p$
is a rational surface and
$\#X_p(\bbF_{\!p}) \equiv 1 \pmod p$
is~automatic. Thus,~let us assume the~contrary.

We~fix an auxiliary prime
number~$l$
that is split
in~$\smash{\bbQ(\sqrt{d})}$
and
let~$p$
be a bad, inert~prime.
For~every
prime~$q$
inert
in~$\smash{\bbQ(\sqrt{d})}$,
choose an absolute Frobenius element
$\smash{\Frob_q \in \Gal(\overline\bbQ/\bbQ)}$.
By~Cebotarev, the elements
$\smash{\sigma^{-1} \Frob_q \sigma \in \Gal(\overline\bbQ/\bbQ)}$,
for~$q$
running through the inert primes and
$\sigma$
through
$\smash{\Gal(\overline\bbQ/\bbQ)}$,
are dense in the coset
$\smash{\Gal(\overline\bbQ/\bbQ) \setminus \Gal(\overline\bbQ/\bbQ(\sqrt{d}))}$,
to which
$\Frob_p$
belongs. The~same is still true when restricting to the
primes~$q$,
at which
$X$
has good~reduction.

For~those, we have the congruence
$\smash{\Tr \Frob_{H^2_\et((X_q)_{\overline\bbF_{\!q}}\!, \bbQ_l)} \equiv 0 \pmod q}$.
In~other words,
$$\smash{\textstyle \Tr \frac1q\!\Frob_{H^2_\et((X_q)_{\overline\bbF_{\!q}}\!, \bbQ_l)} = \Tr \frac1q\!\Frob_{q, H^2_\et(X_{\overline\bbQ}, \bbQ_l)} = \Tr \Frob_{q, H^2_\et(X_{\overline\bbQ}, \bbQ_l(1))}\!}$$
is an integer, necessarily within the range
$[-22,22]$.
As~the condition
$\smash{\Tr \frac1q\varphi \in \bbZ \cap [-22,22]}$
defines a Zariski closed subset
of~$\smash{\GL(H^2_\et(X_{\overline\bbQ}, \bbQ_l(1)))}$,
one has
\begin{equation}
\label{eq_nonord}
\smash{\Tr \Frob_{H^2_\et(X_{\overline\bbQ_p}\!, \bbQ_l)} = \Tr \Frob_{p, H^2_\et(X_{\overline\bbQ}\!, \bbQ_l)} \equiv 0 \pmod p \, ,}
\end{equation}
too, cf.~\cite[Expos\'e~XVI, Corollaire~1.6]{SGA4}

Further,~the eigenvalues of\vspace{.1mm}
$\Frob$
on
$\smash{H^2_\et(X_{\overline\bbQ_p}\!, \bbQ_l)}$
are the same as those on
$\smash{H^2_\et(X_{\overline\bbQ_p}\!, \bbQ_p)}$
\cite[Theorem~3.1]{Oc}. In~addition, a main result of
\mbox{$p$-adic}
Hodge theory \cite[Theorem~III.4.1]{Fa} implies,
as~$X$
is~$K3$,
that not more than one of the eigenvalues of
$\Frob$
on
$\smash{H^2_\et(X_{\overline\bbQ_p}\!, \bbQ_p)}$
may be a
\mbox{$p$-adic}
unit, the others being of strictly positive
\mbox{$p$-adic}~valuation.
Under~these circumstances, the congruence (\ref{eq_nonord}) implies that none of the eigenvalues is a
\mbox{$p$-adic}~unit.

For~comparison with the cohomology
$\smash{H^2_\et((X_p)_{\overline\bbF_{\!p}}\!, \bbQ_l)}$
of the singular fiber, the theory of vanishing cycles~\cite[Expos\'es~I, XIII, and XV]{SGA7} applies, as
$X_p$
has only isolated singularities~\cite[Corollaire~2.9]{Il03}. In~our case, it states that
$\smash{H^2_\et((X_p)_{\overline\bbF_{\!p}}\!, \bbQ_l)}$
naturally injects into
$\smash{H^2_\et(X_{\overline\bbQ_p}\!, \bbQ_l)}$.
In particular, the eigenvalues
of~$\Frob$
on
$\smash{H^2_\et((X_p)_{\overline\bbF_{\!p}}\!, \bbQ_l)}$
form a subset of the 22~eigenvalues
of~$\Frob$
on~$\smash{H^2_\et(X_{\overline\bbQ_p}\!, \bbQ_l)}$.

This~shows that all eigenvalues
on~$\smash{H^2_\et((X_p)_{\overline\bbF_{\!p}}\!, \bbQ_l)}$
are of strictly positive
\mbox{$p$-adic}
valuation. Further,~using the Leray spectral sequence together with the proper base change theorem \cite[Ex\-po\-s\'e~XII, Corollaire 5.2.iii)]{SGA4}, one sees that blow-ups do not affect the transcendental part
$\smash{T_l \subset H^2_\et((X_p)_{\overline\bbF_{\!p}}\!, \bbQ_l)}$.
Hence,
$\smash{\#\widetilde{X}_p(\bbF_{\!p}) \equiv 1 \pmod p}$,
for~$\smash{\widetilde{X}_p}$
the minimal resolution of singularities. The~same is true
for~$X_p$.
}
\eop
\end{lem}

\paragraph*{{\em Counting points on degree-2 K3-surfaces.}}\leavevmode\medskip

\noindent
{\em Structure of our samples.}\hskip3mm
We~consider 
$K3$~surfaces
that are given as desingularizations of the double covers of the projective plane, branched over the union of six~lines. One~reason for choosing this particular family is that it is the one studied before by B.~van Geemen \cite[Ex\-am\-ple~3.4]{Ge}. On~the other hand, this family offers computational advantages,~too.

Our~trial computations with all six lines defined over
$\bbQ$
did not lead to any~success. On~the other hand, six lines defined over an
$S_6$-extension
of~$\bbQ$
and forming a Galois orbit would not be easy to~handle. Our~compromise is as~follows.

The~lines are allowed to form three Galois orbits, each of size~two. Assuming the three
\mbox{$\bbQ$-rational}
points of intersection not to be collinear, we may suppose them without restriction to be
$(1\!:\!0\!:\!0)$,
$(0\!:\!1\!:\!0)$,
and
$(0\!:\!0\!:\!1)$.
The~equation of the surface then takes the~form
$$w^2 = q_1(y,z) q_2(x,z) q_3(x,y) \, .$$
This representation is unique up to action of the monomial group. I.e.,~up to permutation and scaling of the variables.

\begin{algo}[Counting points on one surface]
In~order to determine the number of
\mbox{$\bbF_{\!q}$-rational}
points on one surface, we count the points over the
$q$
affine lines of the form
$(1 \!:\! u \!:\! \star)$
and the affine line
$(0 \!:\! 1 \!:\! \star)$
and sum up these~numbers. Finally,~we add
$1$,
as, on each of our surfaces, there is exactly one point lying
above~$e_3$.
\end{algo}

\begin{rem}[Counting points above one line]
It~is easy to count the number of points above the affine line
$L_{x,y} \colon \Ab^1 \rightarrow \Pb^2$\!,
given by
$t \mapsto (x \!:\! y \!:\! t)$.
Observe~that
$q_3$
is constant on this line. Thus, we get a quadratic twist of an elliptic~curve. The~number of points on it is
$q + \chi(q_3(x,y)) \lambda_{x,y}$,
for
\begin{equation}
\label{eq_ellkurv}
\smash{\lambda_{x,y} := \sum_{t \in \bbF_{\!q}} \chi (q_1(y,t) q_2(x,t))}
\end{equation}
and~$\chi$
the quadratic character
of~$\bbF_{\!q}$.
\end{rem}

\begin{str}[Treating a sample of surfaces]
\label{sample_counting}
Our samples are given by three lists of quadratic forms. One list for
$q_1$,
another for
$q_2$,
and third one
for~$q_3$.
In~the case that we want to count the points on all surfaces, given by the Cartesian product of the three lists, we perform as~follows.

\begin{iii}
\item
For each quadratic form
$q_3$,
compute the values of
$\chi(q_3(1,\star))$
and
$\chi(q_3(0,1))$
and store them in a~table.
\item
Run in an iterated loop over all pairs
$(q_1, q_2)$.
For each pair, do the~following.

$\bullet$
Using~formula (\ref{eq_ellkurv}), compute
$\lambda_{1,\star}$
and~$\lambda_{0,1}$.

$\bullet$
Run~in a loop over all
forms~$q_3$.
Each~time, calculate
$S_{q_1,q_2,q_3} := \sum_\star \chi(q_3(1,\star)) \lambda_{1,\star}$,
using the precomputed~values. The~number of points on the surface, corresponding
to~$(q_1,q_2,q_3)$,
is then
$q^2 + q + 1 + \chi(q_3(0,1)) \lambda_{0,1} + S_{q_1,q_2,q_3}$.
\end{iii}
\end{str}

\begin{rems}
i)
(Complexity and performance).
In~the case that the number of quadratic forms is bigger than
$q$,
the costs of building up the tables are small compared to the final~step. Thus,~the complexity per surface is essentially reduced to
$(q+1)$
table look-ups for the quadratic character and
$(q+1)$
look-ups in the small table, containing the values
$\lambda_{1,\star}$
and~$\lambda_{0,1}$.\smallskip

\noindent
ii)
We are limited by the memory transfer generated by the former table access. We store the quadratic
character in an 8-bit signed integer variable. This doubles the speed compared to a 16-bit variable.
\end{rems}

\looseness-1
\begin{rem}[Detecting real multiplication]
We~used the point counting algorithm, in the version described in~\ref{sample_counting}, within the deterministic Algorithm~\ref{alg_det}, in order to detect
$K3$~surfaces
having real multiplication by a prescribed quadratic number~field. This~allowed us to test more than
$2.2 \cdot 10^7$
surfaces per second on one core of a 3.40\,GHz Intel${}^{
\text{(R)}}$Core${}^{\text{(TM)}}$\mbox{i7-3770}~processor.
The~code was written in plain~{\tt C}.
\end{rem}

\paragraph*{{\em The results.}}

i)
A~run of Algorithm~\ref{alg_stat} over all triples
$(q_1,q_2,q_3)$
of coefficient
height~$\leq \!12$,
using the method described in~\ref{sample_counting} for point counting, discovered the first five surfaces that were likely to have real multiplication
by~$\smash{\bbQ(\sqrt{5})}$.
Observe~that a sample of more than
$10^{11}$~surfaces
was necessary to bring these examples to~light.

Analyzing~the examples, we observed that the product of the discriminants of the three binary quadratic forms was always a perfect~square.\smallskip

\noindent
ii)
We~added this restriction to our search strategy, which massively reduces the number of surfaces to be inspected. Doing so, we could raise the search bound up
to~$80$.
This resulted in more surfaces with probable real multiplication
by~$\smash{\bbQ(\sqrt{5})}$
and one example that was likely to have real multiplication
by~$\smash{\bbQ(\sqrt{2})}$.

From~the results, we observed that the square class of one of the three discriminants always coincided with the discriminant of the field of real~multiplication.\smallskip

\noindent
iii)
This~restriction led to a further reduction of the search~space. At~a final stage, we could raise the search bound
to~$200$
for real multiplication by
$\smash{\bbQ(\sqrt{2})}$,
$\smash{\bbQ(\sqrt{5})}$,
$\smash{\bbQ(\sqrt{13})}$,
and
$\smash{\bbQ(\sqrt{17})}$.
We found many more examples for
$\smash{\bbQ(\sqrt{2})}$
and~$\smash{\bbQ(\sqrt{5})}$,
one example
for~$\smash{\bbQ(\sqrt{13})}$,
but none
for~$\smash{\bbQ(\sqrt{17})}$.

\looseness-1
\begin{rem}
The~final sample for
$\smash{\bbQ(\sqrt{17})}$
consisted of about
$4.18 \!\cdot\! 10^{13}$~surfaces
and required about
$24$~days
of CPU~time. The~computations were executed in parallel on two machines, making use of two cores on each~machine. The~other samples were comparable in~size.
\end{rem}

In the cases of
$\smash{\bbQ(\sqrt{2})}$
and
$\smash{\bbQ(\sqrt{5})}$,
the examples found were sufficient to guess
\mbox{$1$-pa}\-ram\-e\-ter
families. To~summarize, our experiments led us to expect the following two~results. For~the first, we could later devise a proof, the second remains a~conjecture.

\begin{obst}
\label{congr_Qw2_familie}
Let\/~$t \in \bbQ$
be arbitrary and\/
$X^{(2,t)}$
be the\/
$K3$~surface
given~by\/
\begin{eqnarray*}
w^2 & = & \textstyle
 [(\frac18 t^2 \!-\! \frac12 t \!+\! \frac14)y^2 + (t^2 \!-\! 2t \!+\! 2)yz + (t^2 \!-\! 4t \!+\! 2)z^2] \\
 & & \textstyle\hspace{4mm}
 [(\frac18 t^2 \!+\! \frac12 t \!+\! \frac14)x^2 + (t^2 \!+\! 2t \!+\! 2)xz + (t^2 \!+\! 4t \!+\! 2)z^2]
 [2x^2 + (t^2 \!+\! 2)xy + t^2y^2] \, .
\end{eqnarray*}
Then\/
$\#X^{(2,t)}_p(\bbF_{\!p}) \equiv 1 \pmod p$
for every prime\/
$p \equiv 3,5 \pmod 8$.\medskip

\noindent
{\bf Proof.}
{\em
The~case
$p=3$
is elementary.
For~$p \neq 3$,
we shall prove this result below in Theorem~\ref{point_count}, under some additional restrictions
on~$t$.
For~the cases left out there, similar arguments~work. Cf.~Remark~\ref{opencases} for a few~details.
}
\eop
\end{obst}

\begin{obsc}
{\rm i)}
Let\/~$t \in \bbQ$
be arbitrary and\/
$X^{(5,t)}$
be the\/
$K3$~surface
given~by\/
$$\textstyle
w^2 =
 [y^2 + tyz + (\frac5{16} t^2 \!+\! \frac54 t \!+\! \frac54)z^2]
 [x^2 + xz + (\frac1{320} t^2 \!+\! \frac1{16} t \!+\! \frac5{16}) z^2]
 [x^2 + xy + \frac1{20} y^2] \, .
$$
Then\/
$\#X^{(5,t)}_p(\bbF_{\!p}) \equiv 1 \pmod p$
for every prime\/
$p \equiv 2,3 \pmod 5$.\smallskip

\noindent
{\rm ii)}
Let\/~$X^{(13)}$
be the\/
$K3$~surface
given~by\/
$$w^2 = (25y^2 + 26yz + 13z^2)(x^2 + 2xz + 13z^2)(9x^2 + 26xy + 13y^2) \, .$$
Then\/
$\#X^{(13)}_p(\bbF_{\!p}) \equiv 1 \pmod p$
for every prime\/
$p \equiv 2,5,6,7,8,11 \pmod {13}$.
\end{obsc}

\begin{rem}
\label{experiment_evidence}
We~verified the congruences above for all
primes~$p < 1000$.
This~concerns
$X^{(13)}$
as well as the
$X^{(5,t)}$,
for any residue class of
$t$
modulo~$p$.

There~is further evidence, as we computed the characteristic polynomials
of~$\Frob_p$
for
$X^{(13)}$
as well as for
$X^{(5,t)}$
and several exemplary values
of~$t \in \bbQ$,
for the good
primes~$p$
below~$100$.
It~turns out that indeed they all show the very particular behaviour, described in Theorem~\ref{factoriz_rat}. To~be concrete, in each case, either
$\smash{\chi_p^\tr}$
is of degree zero, or
$\smash{\chi_{p^f}^\tr}$
is the square of a quadratic polynomial for a suitable positive
integer~$f$,
or
$\smash{\chi_p^\tr}$
is irreducible of degree four, but splits into two factors conjugate
over~$\bbQ(\sqrt{5})$,
respectively~$\bbQ(\sqrt{13})$.
\end{rem}

\section{The proof for real multiplication in the case of the\/
$\smash{\bbQ(\sqrt{2})}$-family}

\begin{lem}
\label{congr_genuegt}
Let\/~$a,D \in \bbZ$
be such that\/
$\gcd(a,D) = 1$
and\/~$X$
a\/~$K3$~surface
over\/~$\bbQ$.
Suppose~that\/
$\#X_p(\bbF_{\!p}) \equiv 1 \pmod p$
for every good
prime\/~$p \equiv a \pmod D$.
Then\/~$X$
has real or complex~multiplication.\medskip

\noindent
{\bf Proof.}
{\em
For~each prime
$p$,
choose an absolute Frobenius element
$\smash{\Frob_p \in \Gal(\overline\bbQ/\bbQ)}$.
By Cebotarev's density theorem, the elements
$\smash{\sigma^{-1} \Frob_p \sigma \in \Gal(\overline\bbQ/\bbQ)}$,
for the good primes
$p \equiv a \pmod D$
and
$\smash{\sigma \in \Gal(\overline\bbQ/\bbQ)}$,
are topologically dense in the coset of
$\smash{\Gal(\overline\bbQ/\bbQ)}$
modulo
$\smash{\Gal(\overline\bbQ/\bbQ(\zeta_D))}$,
they belong~to. Thus,~there are finitely many elements
$\smash{\sigma_1, \ldots, \sigma_k \in \Gal(\overline\bbQ/\bbQ)}$
such that
$$\{\, \sigma_i \sigma^{-1} \Frob_p \sigma \mid i = 1,\ldots,k, \;p \equiv a \pmod D, \;p \text{ good for } X, \;\sigma \in \Gal(\overline\bbQ/\bbQ) \,\}$$
is dense in
$\smash{\Gal(\overline\bbQ/\bbQ)}$.

Now~choose any prime
$l \not\equiv a \pmod D$,
put
$T_l \subset H^2_\et(X_{\overline\bbQ}, \bbQ_l)$
to be the transcendental part of
\mbox{$l$-adic}
cohomology, and write
$r := \dim T_l$.
Then,~for every good prime
$p \equiv a \pmod D$,
one has
$\smash{\Tr \Frob_{p, T_l} = kp}$,
for
$-22 < -r \leq k \leq r < 22$,
and
$\smash{\det \Frob_{p, T_l} = \pm p^r}$.
Hence,
$$\smash{(\Tr \Frob_{p, T_l})^r = \pm k^r \det \Frob_{p, T_l}} \, ,$$
which defines a Zariski closed subset
$\smash{I \subsetneqq \GO(T_l, \langle.\,,.\rangle)}$,
invariant under conjugation. As
$\smash{\GO(T_l, \langle.\,,.\rangle)}$
is~irreducible, the union
$\sigma_1 I \cup \ldots \cup \sigma_k I$
cannot be the whole~group. Consequently,~the image of
$\Gal(\overline\bbQ/\bbQ) \to \GO(T_l, \langle.\,,.\rangle)$
is not Zariski~dense. In~view of Theorem~\ref{TZ}, this is enough to imply real or complex~multiplication.
}
\eop
\end{lem}

\begin{lem}
\label{jac}
Let\/~$C \colon w^2 = F_4(x,y,l)$
be a family of smooth genus-one curves, parametrized by\/
$l \in B$,
for\/~$B$
an integral scheme in characteristic\/
$\neq \!2$
or~$3$,
$c_4(l)$~and\/~$c_6(l)$
its classical invariants, and\/
$\smash{\Delta := \frac{c_4^3(l) - c_6^2(l)}{1728}}$.
Then,~over the open subscheme\/
$D(\Delta) \subseteq B$,
$$I\colon w^2 = x^3 - 27 c_4(l) x - 54 c_6(l)$$
defines a family of elliptic curves, fiber-wise isomorphic to the relative Jacobian
of\/~$C$.\medskip

\noindent
{\bf Proof.}
{\em
The~existence of the relative Jacobian
$\calJ$
follows from~\cite[Expos\'e~232, Th\'e\-o\-r\`eme 3.1]{FGA}. This~is a family of elliptic~curves.
$I$~is
a family of elliptic curves, too, as
$-16[4(-27c_4(l))^3 + 27(-54c_6(l))^2] = 6^{12} \Delta(l) \neq 0$.

Further,~the generic fiber
$I_\eta$
is isomorphic to the Jacobian of
$C_\eta$~\cite[Proposition 2.3]{Fi}. Thus,~over
$D(\Delta)$,
we have two families of elliptic curves that coincide over the generic point
$\eta \in D(\Delta)$.
The~assertion follows from this, since the moduli stack of elliptic curves is separated~\cite[First main Theorem 5.1.1, together with 2.2.11]{KM}.
}
\eop
\end{lem}

\begin{theo}[The point count]
\label{point_count}
Let\/~$\bbF_{\!q}$
be a finite field of characteristic\/
$\neq \!2,3$
such that\/
$2$
is a non-square
in\/~$\bbF_{\!q}$
and\/
$V$
the singular surface given by\/
$w^2 = q_1(y,z) q_2(x,z) q_3(x,y)$,
for\/
$t \in \bbF_{\!q}$~and
\begin{eqnarray*}
 q_1(y,z) & := & \textstyle (\frac18 t^2 \!-\! \frac12 t \!+\! \frac14)y^2 + (t^2 \!-\! 2t \!+\! 2)yz + (t^2 \!-\! 4t \!+\! 2)z^2 \, , \\
 q_2(x,z) & := & \textstyle (\frac18 t^2 \!+\! \frac12 t \!+\! \frac14)x^2 + (t^2 \!+\! 2t \!+\! 2)xz + (t^2 \!+\! 4t \!+\! 2)z^2 \, , \\
 q_3(x,y) & := & (x+y)(2x+t^2y) \, .
\end{eqnarray*}
Suppose~that\/
$t \neq 0$
and\/~$t^2 \neq -2$.
Then~$\#V(\bbF_{\!q}) = q^2 + q + 1$.\medskip

\noindent
{\bf Proof.}
{\em
We~will prove this result in several~steps.\smallskip

\noindent
{\em First step.}
Preparations.

\noindent
We~will count fiber-wise using the fibration, given
by~$y\!:\!x = l$,
for
$l \in \Pb^1(\bbF_{\!q})$.
This~will yield a result by
$q$
too large, as the point lying over
$(0\!:\!0\!:\!1)$
will be counted
$(q+1)$~times.

The~fiber
$V_l$
is the curve, given by
$w^2 = (1+l) (2+t^2 l) x^2 q_1(lx,z) q_2(x,z)$.
A~partial resolution is provided by
$C_l\colon w^2 = (1+l) (2+t^2 l) q_1(lx,z) q_2(x,z)$,
which defines an elliptic~fibration.

We~claim that
$\sum_l \#C_l(\bbF_{\!q}) = \sum_l \#V_l(\bbF_{\!q})$.
Indeed,~the two fibrations differ only over the line
``$x=0$''.
Since~$V$
ramifies over this line,
$V_x$
has exactly
$(q+1)$~points.
On~the other hand, the curve
$C_x$
is given by
$w^2 = (t^4 \!-\! 12t^2 \!+\! 4) \!\cdot\! (1+l)(2+t^2 l)$.
Here,~the constant
$$t^4 - 12t^2 + 4 = (t^2 - 6)^2 - 32$$
is non-zero, as
$2$
is not a~square.
Thus,~$C_x$
is a double cover of
$\Pb^1\!$,
ramified at
$(-1)$
and~$\smash{(-\frac2{t^2})}$.
But~$\smash{-1 \neq -\frac2{t^2}}$,
since
$2$
is a non-square. In~other words,
$C_x$
is a conic, which has exactly
$(q+1)$~points.\smallskip

\noindent
{\em Second step.}
Singular fibers.

\noindent
There are four singular fibers, at
$l = -1$,
$\smash{-\frac2{t^2}}$,
$0$,
and~$\infty$.
In~fact, for the first two, the coefficient is zero, while, for the others, one of the quadratic forms has a double~zero. We~claim that these are the only singular
\mbox{$\bbF_{\!q}$-rational}
fibers.

To~see this, we first observe that
$q_1$
is of discriminant
$$\textstyle (t^2 - 2t + 2)^2 - 4(\frac18 t^2 - \frac12 t + \frac14)(t^2 - 4t + 2) = (t^2 + 2t + 2)^2 - \frac12 (t^2 + 4t + 2)^2 = \frac12(t^2-2)^2$$
and the same
for~$q_2$.
This~term does not vanish, for any value
of~$t$.
Therefore,~$q_1$
and~$q_2$
always define two lines each, never a double~line. Consequently,~for
$l \neq 0,\infty$,
neither of the two quadratic factors
$q_1(lx,z)$
and~$q_2(x,z)$
may have a double~zero.

To~exclude a common zero, one has to compute the resultant, which turns out to~be
$$\textstyle \frac1{64} (t^4 - 12t^2 + 4)^2
    (l^2 + \frac{-6t^4 + 8t^2 - 24}{\phantom{-}t^4 - 12t^2 + 4}l + 1)
    (l^2 + \frac{-2t^4 - 8t^2 -  8}{\phantom{-}t^4 - 12t^2 + 4}l + 1) \, .$$
Here,~$t^4 - 12t^2 + 4 \neq 0$.
Further,~the quadratic polynomials
in~$l$
are of the discriminants
$\smash{\frac{32(t^2-2)^2(t^2+2)^2}{\;\;\;(t^4 - 12t^2 + 4)^2}}$
and
$\smash{\frac{128t^2(t^2-2)^2}{(t^4 - 12t^2 + 4)^2}}$,
which are non-squares
in~$\bbF_{\!q}$,
because of
$t \neq 0$
and~$t^2 \neq \pm2$.
Thus,~the resultant does not vanish for any value of
$l$,
as long as
$t$
is~admissible.\smallskip

\noindent
{\em Third step.}
Points on the singular fibers.

\noindent
The curves
$C_{-1}$
and~$\smash{C_{\!\!-\!\frac2{t^2}}}$
are part of the ramification locus and therefore degenerate to~lines. They~have
$(q+1)$
points~each.

On~the other hand, the fibers
$C_0$
and~$C_\infty$
are given by
$w^2 = 2(t^2 - 4t + 2)z^2 q_2(x,z)$
and
$w^2 = t^2(t^2 + 4t + 2)z^2 q_1(y,z)$.
Both~are conics with the points
over~$z=0$
unified into a double~point. The~corresponding points on the non-singular
conics~$C_0^\ns$
and
$C_\infty^\ns$
satisfy
$\smash{w^2 = \frac14 (t^4 - 12t^2 + 4)}$,
and
$w^2 = \smash{\frac{t^2}8} (t^4 - 12t^2 + 4)$,
respectively. The~two equations differ by a factor
of~$\smash{\frac{t^2}2}$,
which is a non-square. Hence,~one of the curves
$C_0^\ns$
and
$C_\infty^\ns$
has two points such that
$z=0$,
the other~none. Accordingly,~one of the singular curves
$C_0$
and
$C_\infty$
has
$q$~points,
the other~$(q+2)$.

It~therefore remains to show that
$\sum_{l, C_l\,\text{\scriptsize smooth}} \#C_l(\bbF_{\!q}) = (q-3)(q+1)$.\smallskip

\noindent
{\em Fourth step.}
The classical invariants
$c_4$
and~$c_6$.

\noindent
The~invariants
$c_4$
and~$c_6$
of the family of binary quartic forms defining
$C$
are polynomials in
$l$
and~$t$.
They~may easily be written down, but the formulas become quite~lengthy. The~discriminant~$\Delta$
turns out to be
\begin{align*}
\textstyle \Delta = \frac1{1024} t^{12} (t^2-2)^4 (t^4-12t^2+4)^4
           l^2 (l + \frac2{t^2})^6 (l+1)^6 & \\[-1mm]
\textstyle (l^2 + \frac{-6t^4 + 8t^2 - 24}{\phantom{-}t^4 - 12t^2 + 4}l + 1)^2 &
\textstyle (l^2 + \frac{-2t^4 - 8t^2 -  8}{\phantom{-}t^4 - 12t^2 + 4}l + 1)^2 \, .
\end{align*}
The~arguments given in the second step show that
$\Delta \neq 0$,
except for
$l = -1$,\!
$\smash{-\frac2{t^2}}$,\!
$0$,
and~$\infty$.

By~Lemma~\ref{jac},
$I_l\colon w^2 = x^3 - 27 c_4(l) x - 54 c_6(l)$
is isomorphic to the Jacobian
$\Jac C_l$,
for
$l \neq -1, \smash{-\frac2{t^2}}, 0, \infty$.
This~implies
$\#C_l(\bbF_{\!q}) = \#(\Jac C_l)(\bbF_{\!q}) = \#I_l(\bbF_{\!q})$,
since genus-one curves over finite fields always have~points.

We~have to prove that
$\smash{\sum_{l, C_l\,\text{\scriptsize smooth}} \#I_l(\bbF_{\!q}) = (q-3)(q+1)}$.
I.e.,~that the
$(q-3)$
smooth fibers
of~$I$
have, on average, exactly
$(q+1)$~points.\smallskip\pagebreak[3]

\noindent
{\em Fifth step.}
$l$
versus
$\frac1l$.

\noindent
For~the
\mbox{$j$-invariant}
$\smash{j = \frac{c_4^3}\Delta}$,
one computes that
$\smash{j(\frac1l) = j(l)}$.
More~precisely,
$$\textstyle c_4(\frac1l) = K^2 c_4(l)
\quad\text{and}\quad
c_6(\frac1l) = K^3 c_6(l) \, ,$$
for~$\smash{K := \frac{2l + t^2}{l^4(t^2l + 2)}}$.

In~other words, the elliptic curves
$I_l$
and
$\smash{I_{\frac1l}}$
are geometrically isomorphic to each~other. They~are quadratic twists, according to the extension
$\smash{\bbF_{\!q}(\sqrt{K})/\bbF_{\!q}}$.
Consequently,~if
$\smash{\frac{\,2l + t^2}{t^2l + 2} \in \bbF_{\!q}}$
is a non-square then
$I_l$
and
$\smash{I_{\frac1l}}$
together have exactly
$2(q+1)$~points.\smallskip

\noindent
{\em Sixth step.}
Reparametrization.

\noindent
We~reparametrize according to the M\"obius transformation
$\Pb^1 \to \Pb^1$\!,
$\smash{l \mapsto s := \frac{\,2l + t^2}{t^2l + 2}}$.
This~is not a constant map, for any value
of~$t$.
Indeed,~the determinant of the corresponding
$2\times2$-matrix
is
$4-t^4 = (2-t^2)(2+t^2) \neq 0$.
The~inverse transformation is given by
$\smash{s \mapsto l := \frac{-2s + t^2}{\;t^2\!s - 2}}$.

Write~$I'$
for the fibration, defined by
$\smash{I'_s := I_l}$.
Then~the bad fibers are located at
$\smash{s = -1, \infty, \frac{t^2}2, \frac2{t^2}}$.
The correspondence
$\smash{l \mapsto \frac1l}$
goes over into
$$\textstyle s = \frac{\,2l + t^2}{t^2l + 2} \mapsto \frac{\,\frac2l + t^2}{\frac{t^2}l + 2} = \frac{2 + t^2 l}{t^2 + 2l} = \frac1s \, .$$
Thus,~for
$\smash{s \neq -1, \frac{t^2}2, \frac2{t^2} \in \bbF_{\!q}^*}$
a non-square, the fibers
$I'_s$
and~$\smash{I'_{\frac1s}}$
together have exactly
$2(q+1)$~points.\smallskip

\noindent
{\em Seventh step.}
Pairing the squares~I.

\noindent
It~remains to consider the fibers
for~$s \in \bbF_{\!q}^*$,
$s \neq -1$,
a square and
for~$s=0$.
For~these,
$\smash{I'_s \cong I'_{\frac1s}}$,
except
for~$s = 0$.
There~are
$4n+2$
such fibers, for
$q = 8n+3$
as well as
for~$q = 8n+5$.

It~turns out that
$\smash{j'(s_2) = j'(s_1)}$,
for~$s_1 = a^2$
and
$\smash{s_2 = \frac{(a-1)^2}{(a+1)^2}}$.
More~precisely,
$$\textstyle c'_4(s_2) = F^2 c'_4(s_1)
\quad\text{and}\quad
c'_6(s_2) = F^3 c'_6(s_1) \, ,$$
for
$$\textstyle F := 8\,\frac{(a+1)^2 (a^2 - \frac2{t^2})^4}{(a^2 - \frac{2t^2+4}{t^2-2}a + 1)^4} \, .$$
We~observe here that the denominator never vanishes for
$t \neq 0$,
In~fact, the discriminant of the quadratic polynomial is equal
to~$\smash{\frac{32t^2}{(t^2-2)^2}}$,
which is always a non-square.
As~$8$
is a non-square, we see that
$F$~is
a non-square as long as
$F \neq 0$,
which happens to be true
for~$a \neq -1$.

In~other words,
for~$a \neq -1$,
the elliptic curves
$I'_{s_1}$
and
$\smash{I'_{s_2}}$
are non-trivial quadratic twists of each~other. This~shows that
$\#I'_{s_1}(\bbF_{\!q}) + \#I'_{s_2}(\bbF_{\!q}) = 2(q+1)$.\smallskip

\noindent
{\em Eighth step.}
Pairing the squares~II.

\noindent
In~particular,  we have
$\#I'_1(\bbF_{\!q}) + \#I'_0(\bbF_{\!q}) = 2(q+1)$.
For~the other
$4n$~fibers,
we argue as~follows.
The~group
$V := \bbZ/2\bbZ \times \bbZ/2\bbZ$
operates
on~$\Pb^1(\bbF_{\!q})$
via
$e_1 \!\cdot\! a := -a$
and~$e_2 \!\cdot\! a := \frac1a$.
The~orbits are of size four, except for
$\{0,\infty\}$,
$\{1,-1\}$,
and, possibly,
$\{i,-i\}$.
The~map
$I\colon \Pb^1(\bbF_{\!q}) \to \Pb^1(\bbF_{\!q}), \smash{a \mapsto \frac{a-1}{a+1}}$,
is compatible with the operation
of~$V$
in the sense that
$e_1 \!\cdot\! I(a) = I(e_2 \!\cdot\! a)$
and~$e_2 \!\cdot\! I(a) = I(e_1 \!\cdot\! a)$.

Therefore,~$I$
defines a mapping
$\smash{\overline{I}\colon \Pb^1(\bbF_{\!q})/V \to \Pb^1(\bbF_{\!q})/V}$
from the orbit set to~itself. One~easily sees that
$I(I(a)) = e_1e_2 \!\cdot\! a$.
I.e.,~$\smash{\overline{I}}$
is actually an~involution. Solving~the equations
$\smash{\frac{a-1}{a+1} = \pm a}$
and
$\smash{\frac{a-1}{a+1} = \pm \frac1a}$,
utilizing the fact that
$2$
is a non-square, we find that
$I$
has no fixed points, except for the possible
orbit~$\{i,-i\}$.

Accordingly,~$\smash{J\colon a^2 \mapsto (\frac{a-1}{a+1})^2}$
defines an involution of the squares
in~$\Pb^1(\bbF_{\!q})$
modulo the equivalence relation generated by
$\smash{x \sim \frac1x}$.
The~only possible fixed point
of~$J$
is~$\{-1\}$.
Further,~$J(\{0,\infty\}) = \{1\}$.

As~a consequence, we see that the squares
$x \in \bbF_{\!q}^*$,
different from
$\pm1$,
decompose into sets
$\smash{\{a^2, \frac1{a^2}, (\frac{a-1}{a+1})^2, (\frac{a+1}{a-1})^2\}}$
of exactly four~elements. The~assertion follows immediately from this.
}
\eop
\end{theo}

\begin{rem}
\label{opencases}
If~$t^2 = -2$
then the same result is~true.
For~$t=0$,
however, one has
$\#V(\bbF_{\!q}) = q^2 + 2q + 1$,
while, for
$t=\infty$,
$\#V(\bbF_{\!q}) = q^2 + 1$.
Only~minor mod\-i\-fi\-ca\-tions of the argument are~necessary. The~case
$t^2 = -2$
is actually simpler, as then
$K = -1$
is constant and easily seen to be a non-square. In~each case, there are exactly four singular
\mbox{$\bbF_{\!q}$-rational}
fibers.
\end{rem}

\begin{rems}
i)
Elliptic
$K3$~surfaces
generally have 24 singular~fibers. In~our case,
$I_{-1}$
and~$\smash{I_{\!-\!\frac2{t^2}}}$
are of Kodaira type
$\Ib_0^*$,
thus being of multiplicity~six. The~other six singular fibers,~four of which are defined only
over~$\smash{\bbF_{\!q^2}}$,
are of Kodaira type
$\Ib_2$
and multiplicity~two.\smallskip

\noindent
ii)
The~symmetry under
$\smash{l \leftrightarrow \frac1l}$
is enforced by the~construction. In~fact, consider the double cover
of~$\Pb^2$\!,
branched over the union of the four lines
$z=a_1x$,
$z=a_2x$,
$z=b_1y$,
and~$z=b_2y$.
The~fiber for
$y \!:\! x = l$
has branch points at
$a_1,a_2,b_1l,b_2l$,
which is a quadruple projectively equivalent to
$\smash{a_1,a_2,Kb_1/l,Kb_2/l}$,
for~$\smash{K := \frac{a_1a_2}{b_1b_2}}$.
For~our fibration, independently of the
parameter~$t$,
we have
$\smash{K = \frac{q_2(1,0)}{q_2(0,1)} \!:\! \frac{q_1(1,0)}{q_1(0,1)} = 1}$.

The~twist factor is
$\smash{q_3(1,K/l) / q_3(1,l)}$.
This~expression would be fractional-quadratic, in general, but is fractional-linear in our~case.\smallskip

\noindent
iii)
We~found the second symmetry, which allowed us to pair the squares, by looking at the factorizations of the rational functions
$j(l) - C$.
It~seems to be very specific for the particular fibrations, occurring in the proof of Theorem~\ref{point_count}.
\end{rems}

\begin{theo}[A family of
$K3$~surfaces
with real multiplication]
\label{RM_Qw2_familie}
Let\/~$t \in \bbQ$
be such that\/
$\nu_{17}(t-1) > 0$
and\/
$\nu_{23}(t-1) > 0$.
Then~the\/
$K3$~surface\/
$X^{(2,t)}$
has geometric Picard
rank\/~$16$
and real multiplication
by\/~$\bbQ(\sqrt{2})$.\medskip

\noindent
{\bf Proof.}
{\em
We~proved
$\smash{\#X^{(2,t)}_p(\bbF_{\!p}) \equiv 1 \pmod p}$
for all primes
$p \equiv 3,5 \pmod 8$,
$p > 3$,
in~Theorem~\ref{point_count}. By~Lemma~\ref{congr_genuegt}, this guarantees that
$\smash{X^{(2,t)}}$
has real or complex mul\-ti\-pli\-ca\-tion by a number
field~$E$.

Further,~all the surfaces
$X^{(2,t)}$
considered coincide
modulo~$17$
and
modulo~$23$,
these two primes being~good. Counting~points, one finds
$\smash{\#X^{(2,t)}_{17}(\bbF_{\!17^i}) = 313}$,
$83\,881$,
and~$24\,160\,345$,
as well as
$\smash{\#X^{(2,t)}_{23}(\bbF_{\!23^i}) = 547}$,
$280\,729$,
and~$148\,114\,771$,
for~$i=1,2,3$.
The~characteristic polynomials of
$\Frob_{17}$
and
$\Frob_{23}$
turn out to be
\begin{align*}
\chi^\tr_{17}(Z) = {}& Z^4 + 28Z^3 + \phantom{0}646Z^2 + \phantom{0}8092Z + \phantom{0}83521 \qquad\text{and} \\
\chi^\tr_{23}(Z) = {}& Z^4 + 52Z^3 + 1702Z^2 + 27508Z + 279841 \, ,
\end{align*}
both being~irreducible. In~particular,
$\smash{\rk\Pic(X^{(2,t)}_{\overline\bbF_{\!17}}) = \rk\Pic(X^{(2,t)}_{\overline\bbF_{\!23}}) = 18}$.
Ap\-pli\-ca\-tions of the Artin-Tate formula~\cite[Theorem~6.1]{Mi} show
$$\disc\Pic(X^{(2,t)}_{\overline\bbF_{\!17}}) \in (2 \bmod (\bbQ^*)^2)
\quad\text{and}\quad
\disc\Pic(X^{(2,t)}_{\overline\bbF_{\!23}}) \in (14 \bmod (\bbQ^*)^2) \, .$$
From~this information, one deduces that
$\smash{\rk\Pic(X^{(2,t)}_{\overline\bbQ}) = 16}$
or~$17$.
If~the rank was
$17$
then \mbox{\cite[Theorem~1, together with Remark~2]{Ch1}} shows that
$\smash{\rk\Pic(X^{(2,t)}_{\overline\bbQ}) \leq \rk\Pic(X^{(2,t)}_{\overline\bbF_{\!17}}) - [E:\bbQ]}$,
a contradiction as the right hand side is at
most~$16$.

Our~next assertion is that
$[E:\bbQ] = 2$.
As~$\dim T = 6$,
the potential alternative degrees would be
$3$
or~$6$.
In~the first case,
$E$
is certainly totally~real. In the second case, in view of~\cite[Remark~1.5.3.c)]{Za83},
$E$
must be~CM. In~both cases, there is a totally real, cubic number
field~$E'$,
contained
in~$\End(T)$.

For~$l$
a prime that is inert
in~$E'$,
$T_l$
carries the structure of a vector space over the
field~$E' \!\otimes_\bbQ\! \bbQ_l$.
Further,~there is a
constant~$f$
such that
$\smash{(\Frob_p)^f}$
is an
\mbox{$E' \!\otimes_\bbQ\! \bbQ_l$-linear}
map, for every prime
$p \neq l$.
This,~however, implies that the number of eigenvalues
of~$\smash{(\Frob_p)^f}$,
considered as a
\mbox{$\bbQ_l$-linear}
map, that are roots of unity multiplied
by~$p$,
is a multiple
of~$3$.
The~calculations shown above for
$p=17$
and~$p=23$
clearly disagree with~that.

It~remains to determine the quadratic number
field~$E$
exactly. For~this, an easy computation reveals that the Galois group of
$\chi^\tr_{17}(Z) = Z^4 + 28Z^3 + 646Z^2 + 8092Z + 83521$
is cyclic of order~four. In~particular, variant~i) of Theorem~\ref{factoriz_rat} applies, showing that
$\chi^\tr_{17}$
splits
over~$E$
into two conjugate~factors.
But~$\smash{\bbQ(\sqrt{\disc\chi^\tr_{17}})}$
is the only quadratic subfield of the splitting field
of~$\chi^\tr_{17}$.
A~direct calculation yields, finally,
that~$\disc\chi^\tr_{17} = 2^{29} \!\cdot\! 17^6$.
}
\eop
\end{theo}

\appendix

\section{The analytic approach}

\begin{prop}
\label{Quadsum}
Let\/~$T$
be a\/
\mbox{$\bbQ$-vector}
space of dimension six, equipped with a non-degenerate symmetric, bilinear pairing\/
$\langle.\,,.\rangle \colon T \times T \to \bbQ$
of
discriminant\/
\mbox{$(1 \bmod (\bbQ^*)^2)$}
and\/
$\varphi\colon T \to T$
be a self-adjoint endomorphism such that\/
$\varphi \circ \varphi = [d]$.\smallskip

\noindent
Then\/~$d \in \bbQ$
is a sum of two rational~squares.\medskip

\noindent
{\bf Proof.}
{\em
The proposition is immediate when
$d$
is a~square. Thus,~assume that
$d$
is a non-square. The assumptions
on~$\varphi$
imply that
$\smash{\varphi_{\bbQ(\sqrt{d})}}$
is diagonalizable. For~the eigenvalues
$\smash{\pm\sqrt{d}}$,
the eigenspaces, which we will denote
by~$T_+$
and
$T_-$,
both must be three-di\-men\-sional.
As~$\varphi$
is self-adjoint, they are perpendicular to each~other.

In~particular, the pairings
$\langle.\,,.\rangle |_{T_+}$
and
$\langle.\,,.\rangle |_{T_-}$
are non-degenerate,~too. We~may choose an orthogonal system
$\{x_1, x_2, x_3\} \subset T_+$
such that
$\langle x_i, x_i \rangle =: a_i \neq 0$,
for~$i = 1,2,3$.
Then~the real conjugates
$x'_1, x'_2, x'_3 \in T_-$
also form an orthogonal system, and one has
$\langle x'_i, x'_i \rangle = a'_i \neq 0$.

From~this, one finds an orthogonal decomposition
$T = T_1 \oplus T_2 \oplus T_3$,
defined
over~$\bbQ$,
when~putting
$$T_i := \spann\big(x_i+x_i', \sqrt{d}\, (x_i-x_i')\big) \, .$$
The discriminant
of~$T_i$
is in the class~of
$$\det \!\left(\!\!
\begin{array}{cc}
   a_i+a_i'           & \sqrt{d}\, (a_i-a_i') \\
\sqrt{d}\, (a_i-a_i') &        d(a_i+a_i')
\end{array}
\!\!\right)\!
=
d[(a_i+a_i')^2 - (a_i-a_i')^2] = 4da_ia_i' = 4d N\!_{\bbQ(\!\sqrt{d})\!/\!\bbQ}\!(a_i)
$$
modulo~squares. Consequently,
$\smash{\disc T = ((4d)^3 N_{\bbQ(\!\sqrt{d})/\bbQ}(a_1a_2a_3) \bmod (\bbQ^*)^2)}$.
By our assumption about
$\disc T$,
this implies that
$d$
is a norm
from~$\bbQ(\sqrt{d})$.

As~$(-d)$
is clearly a norm, we conclude that
$(-1)$
must be a norm
from
$\smash{\bbQ(\sqrt{d})}$,
too.
I.e.,~$-1 = a^2 - db^2$
for suitable
$a,b \in \bbQ$.
Therefore,
$d$
is a sum of two~squares.
}
\eop
\end{prop}

\begin{rem}[{cf.~\cite[Example~3.4]{Ge}}]
\label{Ex_vGeemen}
Suppose~$T \cong \bbQ^6$
and that
$\langle.\,,.\rangle$
is the bilinear form defined by the matrix
$\diag(1,1,-1,-1,-1,-1)$.
Then,~for every
$d \in \bbQ$
being a sum of two squares, there exists a self-adjoint endomorphism
$\varphi\colon T \to T$
such that
$\varphi \circ \varphi = [d]$.

Indeed,~decompose
$T$
orthogonally as
$\bbQ^2 \oplus \bbQ^2 \oplus \bbQ^2$
such that, on each summand, the bilinear form is given by either
$\diag(1,1)$
or~$\diag(-1,-1)$.
Then~define
$\varphi$
component-wise by taking the
matrix~$\smash{(\genfrac{}{}{0pt}{1}{u \;\;\;\, v\,}{v\; -u})}$,
for
$d = u^2 + v^2$,
three~times. The~symmetry of the matrix implies that
$\varphi$
is self-adjoint and
$\varphi \circ \varphi = [d]$
is~obvious.
\end{rem}

\begin{theo}
\label{fam_HS}
Let\/~$d \in \bbQ$
be a non-square.

\begin{iii}
\item
If\/~$d$
is not a sum of two squares then there is no
\mbox{weight-$2$}
Hodge structure of dimension six, having a polarization of discriminant\/
\mbox{$(1 \bmod (\bbQ^*)^2)$}
and an endomorphism algebra
containing\/~$\smash{\bbQ(\sqrt{d})}$.
\item
Suppose~that\/
$d$
is a sum of two~squares. Then~there exists a one-di\-men\-sional family of polarized, six-dimensional
\mbox{weight-$2$}
Hodge structures of
$K3$~type,
having the underlying quadratic space
$(\bbQ^6, \diag(1,1,-1,-1,-1,-1))$
and real multiplication
by\/~$\bbQ(\sqrt{d})$.
\end{iii}\medskip

\noindent
{\bf Proof.}
{\em
i)~follows immediately from Proposition~\ref{Quadsum}. Cf.~\cite[Theorem~1.6.a) and Theorem 1.5.1]{Za83}.\smallskip

\noindent
ii)
To~convert
$T := (\bbQ^6, \diag(1,1,-1,-1,-1,-1))$
into a
\mbox{weight-$2$}
Hodge structure of
$K3$~type,
one has to select a one-dimensional isotropic subspace
$H^{2,0} \subset T_\bbC$
such that
$\smash{\overline{H^{2,0}}}$
is not perpendicular
to~$H^{2,0}$.
This~will automatically fix
$\smash{H^{0,2} := \overline{H^{2,0}}}$
and
$\smash{H^{1,1} := (H^{2,0} + \overline{H^{2,0}})^\perp}$.

In~addition, we choose the endomorphism
$\varphi\colon T \to T$
constructed in Remark~\ref{Ex_vGeemen}. By~construction,
$\varphi_\bbC$
commutes with complex conjugation
on~$T_\bbC$.
Furthermore,~as
$\varphi_\bbC$
is self-adjoint and fulfills
$\varphi_\bbC \circ \varphi_\bbC = [d]$,
it respects~orthogonality. Therefore,
$\varphi_\bbC(H^{2,0}) \subseteq H^{2,0}$
alone will be sufficient for
$\varphi$
to cause real~multiplication.

To~ensure this, let us take
$\smash{H^{2,0} \subset T_{\bbC, +}}$.
The~eigenspace
$T_{\bbC, +}$
has a real basis, given by
$\smash{e_i - \frac{u-\sqrt{d}}v e_{i+1}}$,
for~$i = 1, 3, 5$.
In~this basis, the pairing
$\langle.\,,.\rangle |_{T_{\bbC,+}}$
is given by the non-degenerate matrix
$\smash{\diag(1+(\frac{u-\sqrt{d}}v)^2, -1-(\frac{u-\sqrt{d}}v)^2, -1-(\frac{u-\sqrt{d}}v)^2)}$,
which is indefinite. Consequently,~on
$\Pb(T_{\bbC, +}) \cong \Pb^2\!$,
the condition
$\langle x, x \rangle = 0$
defines a
conic~$C$
and, on this conic,
$\langle x, \overline{x} \rangle \neq 0$
is fulfilled on a dense open~subset.
}
\eop
\end{theo}

\begin{rem}
\label{Familie_vGeemen}
Consider~the four-di\-men\-sional family of
$K3$~surfaces
that are given as desingularizations of the double covers of
$\Pb^2\!$,
branched over the union of six~lines. Then
$\rk\Pic(\frakX) \geq 16$
and we are particularly interested in the surfaces for which equality~occurs.
 
In~any case, the pull-back of a general line and the 15 exceptional curves generate a sub-Hodge structure
$P'$
of
dimension~$16$.
The~symmetric, bilinear form
on~$P'$
is given by the matrix
$\diag(2,-2,\ldots,-2)$.
Indeed,~the exceptional curves have self-intersection
number
$(-2)$
\mbox{\cite[Proposition VIII.13.i)]{Be}}.
According to \cite[Ch.\,IV, Theorem~9]{Se70}, there is an isometry
$P' \cong (\bbQ^{16}, \diag(1,-1,\ldots,-1))$.
\end{rem}

\begin{coro}
Let\/~$d \in \bbQ$
be a non-square being the sum of two~squares. Then~there exists a one-dimensional family of\/
$K3$~surfaces
over\/~$\bbC$,
the generic member of which has Picard rank\/
$16$
and real multiplication
by\/~$\smash{\bbQ(\sqrt{d})}$.\medskip

\noindent
{\bf Proof.}
{\em
As~a quadratic space,
$H = H^2(\frakX,\bbQ)$
is the same for all
$K3$~surfaces.
One~has
$H \cong (\bbQ^{22}, \diag(1,1,1,-1,\ldots,-1))$.
By~\cite[Corollary~14.2]{BPV}, cf.~\cite[Ch.\,IX, Theorem~4]{Sh}, there exists a complex-analytic
$K3$~surface
$\frakX$
for every choice of a one-dimensional subspace
$\spann(x) \subset H_\bbC$
fulfilling
$\langle x, x \rangle = 0$
and~$\langle x, \overline{x} \rangle > 0$.

We~choose
$P' \subset H_\bbC$
as in Remark~\ref{Familie_vGeemen}, put
$T' := (P')^\perp$,
and restrict considerations to
subspaces~$\spann(x) \subset T' \subset H_\bbC$.
By~the classification of the quadratic forms
over~$\bbQ$~\cite[Ch.\,IV, \S3]{Se70},
we have
$T' \cong (\bbQ^6, \diag(1,1,-1,-1,-1,-1))$.

Therefore,~Theorem~\ref{fam_HS} guarantees the existence of a one-dimensional family of subspaces
$\spann(x) \subset T'$
such that
$\langle x, x \rangle = 0$
and
$\langle x, \overline{x} \rangle \neq 0$.
The~construction given shows that the first condition actually defines a
conic~$C$
and that
$\langle x, \overline{x} \rangle > 0$
is satisfied on a non-empty open subset
of~$C$.

We~still have to show that, generically,
$\rk \Pic (\frakX) = 16$.
For~this observe, by the Lefschetz theorem on
\mbox{$(1,1)$-classes},
Picard
rank~$16$
is equivalent to
$\bbQ^6 \cap H^{1,1} = 0$.
To~investigate this condition, let
$0 \neq v = (v_1,\ldots,v_6) \in \bbQ^6$
be any vector. The~inclusion
$v \in H^{1,1}$,
for a particular choice
of~$x$,
implies that
$v \in \spann(x)^\perp$.
I.e.,
$$v_1x_1 + v_2x_2 - v_3x_3 - \ldots - v_6x_6 = 0 \, .$$
This~hyperplane meets the
conic~$C$
in at most two~points. Indeed,~the plane
$\Pb(T_{\bbC, +})$
is not contained in any
\mbox{$\bbQ$-rational}
hyperplane, as an inspection of the base vectors given above immediately~shows. In~total, there are only countably many exceptions, for which
$\rk \Pic (\frakX) > 16$.
}
\eop
\end{coro}


\frenchspacing
\begin{thebibliography}{99}
\bibitem{BPV}
Barth, W., Peters, C., and Van de Ven, A.: Compact complex surfaces, {\em Springer,} Berlin, Heidelberg, New York, Tokyo~2004
\bibitem{Be}
Beauville, A.: Surfaces alg\'ebriques complexes,
Ast\'erisque~54, {\em Soci\'et\'e Math\'e\-matique de France,} Paris~1978
\bibitem{BO}
Berthelot, P.\ and Ogus, A.: Notes on crystalline cohomology, {\em Princeton University Press,} Princeton~1978
\bibitem{Ch1}
Charles, F.: On the Picard number of
$K3$~surfaces
over number fields, {\tt arXiv:1111.4117}
\bibitem{Ch2}
Charles, F.: The Tate conjecture for
$K3$~surfaces
over finite fields, {\tt arXiv:1206.4002}
\bibitem{Co}
Cox, D.\,A.: Primes of the form
$x^2+ny^2$.
Fermat, class field theory and complex multiplication, {\em John Wiley \& Sons,} New York~1989
\bibitem{De71}
Deligne, P.: Th\'eorie de Hodge~II, {\em Publ.\ Math.\ IHES\/} {\bf 40}\br(1971)\brr5--57
\bibitem{De74}
Deligne, P.: La conjecture de Weil~I, {\em Publ.\ Math.\ IHES\/}
{\bf 43}\br(1974)\brr273--307
\bibitem{DM}
Dixon, J.\,D.\ and Mortimer, B.: Permutation groups, Graduate Texts in Mathematics~163, {\em Springer,} New York~1996
\bibitem{EK}
Elkies, N.\ and Kumar, A.:
$K3$~surfaces
and equations for Hilbert modular surfaces, {\tt arXiv:\discretionary{}{}{}1209.3527}
\bibitem{EJ08}
Elsenhans, A.-S.\ and Jahnel, J.:
$K3$~surfaces
of Picard rank one and degree two, in: Algorithmic number theory (ANTS~8), Lecture Notes in Computer Science~5011, {\em Springer,} Berlin~2008, 212--225
\bibitem{EJ10}
Elsenhans, A.-S.\ and Jahnel, J.: On Weil polynomials of
$K3$~surfaces,
in: Algorithmic Number Theory (ANTS~9), Lecture Notes in Computer Science~6197, {\em Springer,} Berlin~2010, 126--141
\bibitem{EJ12}
Elsenhans, A.-S.\ and Jahnel, J.: Kummer surfaces and the computation of the Picard group, {\em LMS Journal of Computation and Mathematics\/} {\bf 15}\br(2012)\brr84--100
\bibitem{EJ13}
Elsenhans, A.-S.\ and Jahnel, J.: On the computation of the Picard group for certain singular quartic surfaces, {\em Mathematica Slovaca\/} {\bf 63}\br(2013)\brr215--228
\bibitem{Fa}
Faltings, G.:
$p$-adic
Hodge theory, {\em J.\ Amer.\ Math.\ Soc.} {\bf 1}\br(1988)\brr255--299
\bibitem{Fi}
Fisher, T.: The invariants of a genus one curve, {\em Proc.\ Lond.\ Math.\ Soc.} {\bf 97}\br(2008)\brr753--782
\bibitem{FKRS}
Fit\'e, F., Kedlaya, K., Rotger, V., and Sutherland, A.\,V.: Sato-Tate distributions and Galois endomorphism modules in genus
$2$,
{\em Compos.\ Math.} {\bf 148}\br(2012)\brr1390--1442
\bibitem{FGA}
Grothendieck, A.: Fondements de la G\'eom\'etrie Alg\'ebrique (FGA), {\em S\'eminaire Bourbaki\/} 149, 182, 190, 195, 212, 221, 232, 236, Paris~1957--62
\bibitem{Ge}
van Geemen, B.: Real multiplication on
$K3$~surfaces
and Kuga-Satake varieties, {\em Michigan Math.\ J.} {\bf 56}\br(2008)\brr375--399
\bibitem{Ha}
Hall Jr., M.: Combinatorial theory, {\em Blaisdell Publishing Co.,} Waltham, Toronto, London~1967
\bibitem{HSBT}
Harris, M., Shepherd-Barron, N., and Taylor, R.: A family of Calabi-Yau varieties and potential au\-to\-mor\-phy, {\em Ann.\ of Math.} {\bf 171}\br(2010)\brr779--813
\bibitem{Hu}
Humbert, G.: Sur les fonctions ab\'eliennes singuli\`eres, {\em J.\ Math.\ Pures et Appl.}
5$^{e}$~s\'erie,
{\bf 5}\br(1899)\brr233--350
\bibitem{Il91}
Illusie, L.: Crystalline cohomology, in: Motives (Seattle~1991), Proc.\ Sympos.\ Pure Math.\ 55-1, {\em AMS,} Providence~1994, 43--70
\bibitem{Il03}
Illusie, L.: Perversit\'e et variation, {\em Manuscripta Math.} {\bf 112}\br(2003)\brr271--295
\bibitem{IR}
Illusie, L.\ and Raynaud, M.: Les suites spectrales associ\'ees au complexe de de Rham-Witt, {\em Publ.\ Math.\ IHES\/} {\bf 57}\br(1983)\brr73--212
\bibitem{KM}
Katz, N.\,M.\ and Mazur, B.: Arithmetic moduli of elliptic curves,
Annals of Mathematics Studies~108, {\em Princeton University Press,} Princeton~1985
\bibitem{KS}
Kedlaya, K.\,S.\ and Sutherland, A.V.: Hyperelliptic curves,
\mbox{$L$-polynomials,}
and random matrices, in: Arithmetic, geometry, cryptography and coding theory, Contemp.\ Math.~487, {\em AMS,} Providence~2009, 119--162
\bibitem{LP}
Larsen, M.\ and Pink, R.: On
\mbox{$l$-independence}
of algebraic monodromy groups in compatible systems of representations,
{\em Invent.\ Math.} {\bf 107}\br(1992)\brr603--636
\bibitem{LMS}
Lieblich, M., Maulik, D., and Snowden, A.: Finiteness of
$K3$~surfaces
and the Tate conjecture, {\tt arXiv:\discretionary{}{}{}1107.1221}
\bibitem{Li}
Liedtke, C.: Lectures on supersingular
$K3$~surfaces
and the crystalline Torelli theorem, {\tt arXiv:1403.2538}
\bibitem{vL05}
van Luijk, R.: Rational points
on~$K3$~surfaces,
{\em Ph.D.\ thesis,} Berkeley~2005
\bibitem{vL07}
van Luijk, R.:
$K3$~surfaces with Picard number one and infinitely many rational points,
{\em Algebra \& Number Theory\/} {\bf 1}\br(2007)\brr1--15
\bibitem{Ma}
Mazur, B.: Frobenius and the Hodge filtration (estimates), {\em Ann.\ of Math.} {\bf 98}\br(1973)\brr58--95
\bibitem{Mi}
Milne, J.\,S.: On a conjecture of Artin and Tate, {\em Ann.\ of Math.} {\bf 102}\br(1975)\brr517--533
\bibitem{Ne}
Neukirch, J.: Algebraic number theory, Grundlehren der Mathematischen Wis\-sen\-schaf\-ten~322, {\em Springer,} Berlin~1999
\bibitem{Oc}
Ochiai, T.:
\mbox{$l$-independence}
of the trace of monodromy, {\em Math.\ Ann.} {\bf 315}\br(1999)\brr321--340
\bibitem{Pe}
Pera, K.\,M.: The Tate conjecture for
$K3$~surfaces
in odd characteristic, {\tt arXiv:1301.6326}
\bibitem{Se70}
Serre, J.-P.: Cours d'arithm\'etique, {\em Presses Universitaires de France,} Paris~1970
\bibitem{SGA4}
Artin, M., Grothendieck, A.\ et Verdier, J.-L.\ (avec la collaboration de Deligne, P.\ et Saint-Donat, B.): Th\'eorie des Topos et Cohomologie \smash{\'Etale} des Sch\'emas, S\'eminaire de G\'eom\'etrie Alg\'ebrique du Bois Marie 1963--1964 (SGA\,4), Lecture Notes in Math.~269, 270, 305, {\em Springer,} Berlin, Heidelberg, New York~1972--1973
\bibitem{SGA5}
Grothendieck, A.\ (avec la collaboration de Bucur, I., Houzel, C., Illusie, L.\ et Serre, J.-P.): Cohomologie
\mbox{$l$-adique}
et
Fonctions~$L$,
S\'eminaire de G\'eom\'etrie Alg\'ebrique du Bois Marie 1965--1966 (SGA\,5), Lecture Notes in Math.~589, {\em Springer,} Berlin, Heidelberg, New York~1977
\bibitem{SGA7}
Deligne, P.\ and Katz, N.: Groupes de Monodromie en G\'eom\'etrie Alg\'ebrique, S\'e\-mi\-naire de G\'eom\'etrie Alg\'ebrique du Bois Marie 1967--1969 (SGA\,7), Lecture Notes in Math.~288, 340, {\em Springer,} Berlin, Heidelberg, New York~1973
\bibitem{Sh}
The members of the seminar of I.\,R.~\v{S}afarevi\v{c}: Algebraic surfaces, {\em AMS,} Providence 1965
\bibitem{Si}
Silverman, J.\,H.: Advanced topics in the arithmetic of elliptic curves,
Graduate Texts in Mathematics~151, {\em Springer,} New York~1994
\bibitem{Ta90}
Tankeev, S.\,G.: Surfaces of
$K3$~type
over number fields and the Mumford-Tate conjecture (Russian), {\em Izv.\ Akad.\ Nauk SSSR Ser.\ Mat.} {\bf 54}\br(1990)\brr846--861
\bibitem{Ta95}
Tankeev, S.\,G.: Surfaces of
$K3$~type
over number fields and the Mumford-Tate conjecture~II (Russian), {\em Izv.\ Ross.\ Akad.\ Nauk Ser.\ Mat.} {\bf 59}\br(1995)\brr179--206
\bibitem{We}
Weber, H.: Lehrbuch der Algebra, 2.~Auflage, 3.~Band: Elliptische Funktionen und algebraische Zahlen, {\em Friedr.\ Vieweg \& Sohn,} Braunschweig~1908
\bibitem{Za83}
Zarhin, Yu.\,G.: Hodge groups of
$K3$~surfaces,
{\em J.\ Reine Angew.\ Math.} {\bf 341}\br(1983)\brr193--220
\bibitem{Za93}
Zarhin,~Yu.\,G.: Transcendental cycles on ordinary
$K3$~surfaces
over finite fields, {\em Duke Math.\ J.} {\bf 72}\br(1993)\brr65-83
\end{thebibliography}
\end{document}